\documentclass[onecolumn,draftcls,10pt]{IEEEtran}
\usepackage{cite}
\usepackage{color}
\usepackage{color, colortbl}
	
\definecolor{LightCyan}{rgb}{0.88,1,1}

\usepackage[linesnumbered,ruled,vlined]{algorithm2e}
\DontPrintSemicolon

\usepackage[pdftex]{graphicx}
\graphicspath{{fig/}{jpeg/}}
\usepackage[cmex10]{amsmath}
\usepackage{amssymb}
\usepackage{multirow}
\usepackage{pifont}% http://ctan.org/pkg/pifont
\newcommand{\cmark}{\ding{51}}%
\newcommand{\xmark}{\ding{55}}%
\usepackage{tikz}
\usetikzlibrary{shapes,arrows,matrix}
\usetikzlibrary{decorations.pathmorphing} % noisy shapes
\usetikzlibrary{fit}					% fitting shapes to coordinates
\usetikzlibrary{backgrounds}	% drawing the background after the foreground
\usepackage{verbatim}
\usepackage{algpseudocode}
\usepackage{amsmath,amssymb,lipsum}
\usepackage{hyperref}
\usepackage{comment}
\usepackage{graphicx}
\usepackage[font=small]{caption}
\usepackage{subcaption}
\usepackage{mathtools,enumerate}
\input{mysymbol.sty}
\usepackage{needspace}

\usepackage{theorem}
\usepackage{bbm}
\usepackage[outdir=./images/]{epstopdf}

\newtheorem{thm}{Theorem}
\newtheorem{lemma}{Lemma}

\newtheorem{corollary}{Corollary}

\theoremstyle{definition}
\newtheorem{assumption}{Assumption}

\def \QED {$\blacksquare$}
\def \w {{\mathbf{w}}}
\def \x {{\mathbf{x}}}

\def \y {{\mathbf{y}}}

\def \eps {\epsilon}
\def \mbE {{\mathbb{E}}}
\def \mbR {{\mathbb{R}}}
\def \bbx {{\mathbf{x}}}
\def \bby {{\mathbf{y}}}
\def \ccalX {{\mathcal{X}}}

\def \L {{\hat{\cL}}}
\def \bbs {{\mathbf{s}}}
\def \bbe {{\mathbf{e}}}

\def \bbd {{\mathbf{d}}}

\def \bbtheta {{\boldsymbol{\theta}}}
\def \bblambda {{\boldsymbol{\lambda}}}
\def \thh {{\theta}}

\def \th {{\theta}}

\def \d {{\bbd}}
\def \cX {{\ccalX}}

\def \M {{\mathbf{M}}}
\def \X {{\mathbf{X}}}
\def \V {{\mathbf{V}}}
\def \Y {{\mathbf{Y}}}
\def \W {{\mathbf{W}}}
\def \Eb {{\mathbf{E}}}
\def \L {{\mathbf{L}}}

\def \n {\nabla}
\def \s {{\bbs}}

\def \e {{\bbe}}

\def \hht {h_{\xi_t}}

\def \ft {f_{\theta_t}}

 \providecommand{\norm}[1]{\left\|#1\right\|}
 \providecommand{\ip}[1]{\boldsymbol{\langle}#1\boldsymbol{\rangle}}

{\tiny }
\title{\vspace{-0cm}{Projection-Free Stochastic Bi-level Optimization}}

\author{Zeeshan~Akhtar, Amrit Singh Bedi, Srujan Teja Thomdapu and Ketan~Rajawat% <-this % stops a space
	 \thanks{
	        Z. Akhtar, S. T. Thomdapu, and K. Rajawat are with the Department of Electrical Engineering,
	        Indian Institute of Technology Kanpur, Kanpur 208016, India (e-mail: zeeshan@iitk.ac.in, srujant@iitk.ac.in,
	        ketan@iitk.ac.in). A. S. Bedi is with the Institute of Systems Research, University of Maryland, College Park, MD, USA. }\vspace{-0mm}}
\begin{document}
\maketitle

\begin{abstract}
	Bi-level optimization, where the objective function depends on the solution of an inner optimization problem, provides a flexible framework for solving a rich class of problems such as hyper-parameter optimization and model-agnostic meta learning. This work puts forth the first \textbf{S}tochastic \textbf{B}i-level \textbf{F}rank-\textbf{W}olfe (SBFW) algorithm to solve the stochastic bi-level optimization problems in a projection-free manner. We utilize a momentum-based gradient tracker that results in a sample complexity of $\mathcal{O}(\eps^{-3})$ for convex objectives and $\mathcal{O}(\eps^{-4})$ for non-convex objectives. The stochastic compositional optimization problems, which are a special case of bi-level optimization problems entailing the minimization of a composition of two expected-value functions, are also considered within the same rubric. The proposed \textbf{S}tochastic \textbf{C}ompositional \textbf{F}rank-\textbf{W}olfe (SCFW) algorithm is shown to achieve a sample complexity of $\mathcal{O}(\eps^{-2})$ for convex objectives and $\mathcal{O}(\eps^{-3})$ for non-convex objectives, at par with the state-of-the-art rates of projection-free algorithms for single-level problems. The usefulness and flexibility of SBFW and SCFW algorithms is demonstrated via extensive numerical tests. We show that SBFW outperforms the state-of-the-art methods for the problem of matrix completion with denoising, and achieve improvements of up to $82\%$ in terms of the wall-clock time required to achieve a particular level of accuracy. Furthermore, we demonstrate the improved performance of SCFW over competing projection-free variants on the policy evaluation problem in reinforcement learning. 
\end{abstract}
%%%%%%%%%%%%%%%%%%%%%%
%%%%%%%%%%%%%%%%%%%%%%
%%%%%%%%%%%%%%%%%%%%%%
%%%Introduction%%%%%%%

\section{Introduction}\label{sec:intro}
We consider the two-level hierarchical optimization problem
\begin{align}\label{main_problem_2}
(\mathcal{P}_1) \ \ \ \ \ &	\min_{\x\in \cX \subset \mathbb{R}^m} Q(\bbx):=F(\bbx, \y^\star(\x)) &(\text{outer})
\\
\ \ \ \ \ \ \ \ \ &\ \ \ \ \ \bby^*(\bbx) \in \arg\min_{\bby}  G(\bby, \bbx)  &(\text{inner}). \label{second-level}
\end{align}
Here,  the outer problem involves minimizing the objective function $Q(\bbx)$ with respect to $\bbx$ over the convex compact constraint set $\mathcal{X}\subset \mathbb{R}^m$.  The objective function is of the form $Q(\bbx):=F(\bbx, \y^\star(\x))$, where $\y^\star(\x)$ is a solution of the inner optimization problem, which for a given $\bbx$, entails minimizing the strongly convex function $G(\bby,\bbx)$ with respect to optimization variable $\bby$. Observe that for bilevel problems of type $\mathcal{P}_1$, the inner and outer problems are inter-dependent and cannot be solved in isolation. Yet, these problems arise in a number of areas, such as meta-learning \cite{rajeswaran2019meta}, continual learning \cite{borsos2020coresets}, reinforcement learning \cite{zhang2020bi}, and hyper-parameter optimization \cite{grazzi2020iteration,franceschi2018bilevel}. Of particular interest are the large-scale or stochastic settings, where the functions $F$ and $G$ are expectations of random functions with unknown distributions, and are accessible only through their samples. First-order stochastic approximation algorithms have been recently proposed \cite{dempe2020bilevel,ghadimi2018approximation,yang2021provably, chen2021single}. The main idea behind these algorithms is to run one gradient descent step in order to solve the inner optimization problem, and subsequently, utilize the updated variable $\bby$ to run another gradient descent step on the outer minimization problem.

\begin{table*}[h]
	\centering\scalebox{0.8}{
		\begin{tabular}{|c|c|c|c|c|c|}
			\hline
			\multirow{1}{*}{{{Reference}}} & \multirow{1}{*}{{Objective}} & \multirow{1}{*}{{Projection Free} } & \multirow{1}{*}{{Problem Type}} & \multirow{1}{*}{{SFO Complexity} {(Outer)}} & \multirow{1}{*}{{SFO Complexity} {(Inner)}} \\\hline
			SFW\cite{mokhtari2020stochastic}       & Convex, Non-Convex                                           & \cmark                                                                     & Single-Level                                                  & $\mathcal{O}(\eps^{-3})$, $\mathcal{O}(\eps^{-4})$                                  & -                                             \\ \hline
			ORGFW  \cite{xie2020efficient}    & Convex                                               & \cmark                                                                     & Single-Level                                                  & $\mathcal{O}(\eps^{-2})$                                    & -                                         \\ \hline
			%SFW \cite{reddi2016stochastic} & Non-Convex                                             & \cmark                                                                      & Single-Level                                                      & $\mathcal{O}(\eps^{-4})$                                  & -                                             \\ \hline SFW$^{++}$\cite{zhang2020one}       & Convex, Non-Convex                                             & \cmark                                                                     & Single-Level                                                & $\mathcal{O}(\eps^{-2})$, $\mathcal{O}(\eps^{-3})$                                   & - \\ \hline 
			SFW$^{++}$\cite{zhang2020one}       & Non-Convex                                           & \cmark                                                                     & Single-Level                                                  & $\mathcal{O}(\eps^{-3})$                                    & -                                              \\ \hline
			SCGD   \cite{wang2017stochastic}                   & Convex, Non-Convex                                               & \xmark                                                                      & Compositional                                                  & $\mathcal{O}(\eps^{-4})$, $\mathcal{O}(\eps^{-4})$                                 & --                                             \\ \hline
			%	SCGD \cite{wang2017stochastic}                   & Non-Convex                                               & \xmark                                                                      & Compositional                                                      & $\mathcal{O}(\eps^{-4})$                                  & --                                           \\ \hline
			ASC  \cite{wang2017accelerating}                                              & Convex, Non-Convex                                                & \xmark                                                                      & Compositional                                                           & $\mathcal{O}(\eps^{-2})$, $\mathcal{O}(\eps^{-2.25})$                                & - , -                                                                \\ \hline
			%		ASC \cite{wang2017accelerating}                                                  & Non-Convex                                               & \xmark                                                                      & Compositional                                                            & $\mathcal{O}(\eps^{-2.25})$                                   & --                                                                \\ \hline
			NASA \cite{ghadimi2020single}, SCSC \cite{chen2020solving}                                                 & non-convex                                               & \xmark                                                                      & Compositional                                                            & $\mathcal{O}(\eps^{-2})$                                   & --                                                                 \\ \hline
			\rowcolor{LightCyan}
			\textbf{SCFW}   (This work)                      & Convex, Non-Convex                                             & \cmark                                                                      & Compositional                                                & $\mathcal{O}(\eps^{-2})$, $\mathcal{O}(\eps^{-3})$                                   & 
			- , -                                         \\ \hline
			%	SCSC \cite{chen2020solving}                                                 & non-convex                                               & \xmark                                                                      & Compositional                                                            & $\mathcal{O}(\eps^{-2})$                                   & --                                                                 \\ \hline
			BSA \cite{ghadimi2018approximation}                    &  non-convex                                              & \xmark                                                                      & Bi-Level                                                    & $\mathcal{O}(\eps^{-2})$                                   & $\mathcal{O}(\eps^{-3})$                                          \\ \hline	stocBiO \cite{ji2020bilevel}, STABLE \cite{chen2021single}, MSTSA \cite{khanduri2021momentum}                   &  non-convex                                              & \xmark                                                                      & Bi-Level                                                    & $\mathcal{O}(\eps^{-2})$                                   & $\mathcal{O}(\eps^{-2})$                                          \\ \hline
			TTSA \cite{hong2020two}                        &  non-convex                                              & \xmark                                                                      & Bi-Level                                                    & $\mathcal{O}(\eps^{-2.5})$                                   & $\mathcal{O}(\eps^{-2.5})$                                          \\ \hline
			%	STABLE \cite{chen2021single}                       &  non-convex                                              & \cmark                                                                      & Bi-Level                                                    & $\mathcal{O}(\eps^{-2})$                                   & $\mathcal{O}(\eps^{-2})$                                          \\ \hline
			SUSTAIN \cite{khanduri2021near}                      &  non-convex                                              & \xmark                                                                      & Bi-Level                                                    & $\mathcal{O}(\eps^{-1.5})$                                   & $\mathcal{O}(\eps^{-1.5})$                                          \\ \hline
			%\textbf{SCFW}                                                      & Non-Convex                                               & \cmark                                                                      & Compositional                                                            & $\mathcal{O}(\eps^{-3})$                                   & -                                                    
			%\\ \hline
			\rowcolor{LightCyan}
			\textbf{SBFW} (This work)                        & Convex, Non-convex                                               & \cmark                                                                      & Bi-Level                                                    & $\mathcal{O}(\eps^{-3})$, $\mathcal{O}(\eps^{-4})$                                   & $\mathcal{O}(\eps^{-1.5})$, $\mathcal{O}(\eps^{-2})$                                          
			%	\\
			%	\textbf{SBFW}                                                      & Non-Convex                                              & \cmark                                                                      & Bi-Level                                                             & $\mathcal{O}(\eps^{-4})$                                    & $\mathcal{O}(\eps^{-2})$                                    
			\\ \hline
	\end{tabular}}
	\caption{ This table compares the different algorithms available in the literature to solve single-level, bi-level, and compositional problems. We present the stochastic first-order complexity (SFO) for the outer as well as inner optimization problems, which denotes the number of minimum queries to the first-order oracle made by the corresponding algorithm to achieve an $\epsilon$ optimal stationary point. Note that the proposed algorithms are able to achieve the best possible SFO complexities as compared to the existing projection-free algorithms.   }
	\label{table1}
	\vspace{-3mm}
\end{table*}
In some works, such as \cite{khanduri2021near,chen2021tighter}, the constraint set $\mathcal{X}$ in the outer optimization problem is taken to be $\mathcal{X}=\mathbb{R}^m$, resulting in a simpler unconstrained outer optimization problem. However, in applications such as meta-learning \cite{rajeswaran2019meta}, personalized federated learning \cite{fallah2020personalized}, and corsets \cite{borsos2020coresets}, the constraint set $\mathcal{X}$ is a strict subset of $\mathcal{X}\subset\mathbb{R}^m$. The standard approach to dealing with such constraint sets is to project the updates of the outer optimization problem onto $\mathcal{X}$ at every iteration. Though popular and widely used, the projected gradient approaches may not necessarily be practical, for instance, in cases where the projection sub-problem is too expensive to be solved at every iteration. The difficulties surrounding projection-based methods have motivated the development of projection-free algorithms \cite{jaggi2013revisiting}, that make use of the Frank-Wolfe (FW) updates \cite{frank1956algorithm}. These FW-based algorithms only require solving a linear program over $\mathcal{X}$, which could be significantly cheaper than solving a non-linear projection problem, as in the case of $\ell_1$-norm or nuclear norm ball constraints. 

Projection-free algorithms for single-level stochastic optimization algorithms are well-known and state-of-the-art algorithms achieve a sample complexity of $\mathcal{O}(\epsilon^{-2})$ \cite{xie2020efficient,zhang2020one,akhtar2022zeroth}. These algorithms rely on a recursive gradient tracking approach that allows the samples to be processed sequentially and achieves variance reduction without the use of checkpoints or large batches. Motivated from these developments, we ask the following question:

``\emph{Is it possible to develop efficient projection-free algorithms for bi-level stochastic optimization problems?}"

This work puts forth the \textbf{S}tochastic \textbf{B}i-level \textbf{F}rank-\textbf{W}olfe (SBFW) algorithm, which is the first projection-free algorithm for bi-level problems. Further, we also focus on the special class of problems called stochastic compositional problems where the inner optimization problem in \eqref{main_problem_2} is solvable in closed-form  \cite{chen2021tighter,khanduri2021near}. Although the results developed for bi-level problems are applicable to compositional problems, we provide improved convergence rates for the stochastic compositional problems using a modified analysis (cf. Sec. \ref{Sec:SCGD}).  
It is important to note that, while bilevel optimization problems are more challenging  to solve, they offer additional flexibility in formulation, and tend to outperform their single-level counterparts. Of particular interests are sparse bilevel problems, such as sparse MAML (\cite{huang2021enhanced,gai2019sparse}), which involve $\ell_1$ or nuclear norm constraints. Other important applications where projection-free bilevel optimization algorithms apply include image reconstruction \cite{crockett2021bilevel} and neural network architecture search \cite{huang2021enhanced}.
Before proceeding, we
provide an practical example and discuss various challenges
associated with solving it efficiently.

\subsection{Motivating Example}
\noindent{\textbf{Matrix Completion with Denoising}}: Let us consider the matrix completion problem with the goal of recovering the missing entries from incomplete and noisy observation of a small and random subset of its entries. In general, for problems without noise, the data matrix can be modeled as a low-rank matrix motivating the use of the nuclear norm constraint \cite{jain2013low}. Low-rank matrix completion problem arises in a wide range of applications such as  image processing \cite{chen2019low}, multi-task learning \cite{evgeniou2007multi},  and collaborative filtering \cite{li2009transfer}. However, under noisy observations, directly solving the matrix completion problem with just the nuclear norm constraints  can result in suboptimal performance \cite{yokota2017simultaneous,mcrae2021low}. Further, noise is present in many vision applications \cite{ji2010robust,keshavan2009matrix}, and using only low-rank priors is not sufficient for recovery of the underlying matrix. A common approach to tackle the noise is to apply a denoising algorithm as a pre-processing step. In general however, it is necessary to apply some heuristics, since denoising algorithms requires access to the full matrix, which is not available in the pre-processing stage. Denoising is also impractical in online settings, where a random subset of the entries of the matrix are observed at every iteration. The bilevel optimization framework provides a way out, allowing the incorporation of denoisining step withing the inner-level subproblem. Mathematically, the  bi-level matrix completion with denoising problem can be written as
\begin{align}\label{mat_comp}
\min_{\norm{\X}_{*}\leq \alpha}&\; \frac{1}{|\Omega_1|}\sum_{(i,j)\in \Omega_1}(\X_{i,j}-\Y_{i,j})^2, \\
&\hspace{-8mm}\text{s. t.}\quad \Y\in\arg \min_{\V} \big\{ \frac{1}{|\Omega_2|}\sum_{(i,j)\in \Omega_2}(\V_{i,j}-\M_{i,j})^2+\lambda_1 \norm{\V}_1 + \lambda_2 \norm{\X-\V}_F^2\big\},\nonumber
\end{align}
%
%\begin{align}\label{mat_comp}
%\min_{\norm{\X}_{*}\leq \alpha}&\mbE\norm{\mathcal{P}_{\Omega}(\X)- \mathcal{P}_{\Omega}(\Y)}^2_F, \\
%&\hspace{-8mm}\text{s. t.}\quad \Y\in\arg \min_{\Y} \big\{ ({1}/{2})\mathbb{E }\norm{\mathcal{P}_{\Omega}(\Y)- \mathcal{P}_{\Omega}(\M)}_F^2\nonumber
%\\
%&\hspace{2.5cm}+\lambda_1 \norm{\Y}_1 + \lambda_2 \norm{\X-\Y}_F^2\big\},\nonumber
%\end{align}
%
where $\M\in \mathbb{R}^{n\times m}$ is the given incomplete noisy matrix, $\norm{\V}_1:=\sum_{i,j}|\V_{i,j}|$ is the sum-absolute-value ($\ell_1$) norm, and $\lambda_1$ and $\lambda_2$ are regularization parameters.   Note that the regularization over the discrepancy between $\X$ and denoised matrix $\Y$ results in bilevel formulation \eqref{mat_comp}. A similar technique in deterministic settings is utilized in various other applications in machine learning and signal processing problems \cite{mccann2020supervised,crockett2021motivating}. 
%The expectation in \eqref{mat_comp} is with respect to the noisy random observations of set $\Omega$.
%
Observe that \eqref{mat_comp} is in the form of bilevel formulation \eqref{main_problem_2}; however, when the entries are revealed in the form of randomly selected subsets $\Omega_1^t \subset \Omega_1$ and $\Omega_2^t \subset \Omega_2$  at every iteration, it becomes stochastic in nature (cf. Sec. \ref{sec:Pre}). The main challenge here is due to the nuclear norm constraint, which makes it quite computationally expensive (sometimes even impractical) to solve \eqref{mat_comp} using projection-based bilevel algorithms. In Sec. \ref{sec:application}, we will show experimentally that the proposed algorithm SBFW is best suited to address such challenges in large-scale bi-level stochastic optimization problems. 

\subsection{Related Work}
We review some of the related work in the context of
bi-level optimization, compositional optimization, and projection-free algorithms.

\noindent{\textbf{Bi-level optimization}} has had a long history, with the earliest applications in economic game theory \cite{von1952theory}. Bi-level optimization has recently received great attention from the machine learning community due to the number of applications in the area  \cite{bennett2008bilevel}. A series of works that proposed to solve the problem of the form \eqref{main_problem_2} has appeared recently \cite{ghadimi2018approximation,yang2021provably,  hong2020two, chen2021single, khanduri2021near, huang2021biadam, huang2021enhanced}. Of these, the seminal works in \cite{ghadimi2018approximation,yang2021provably} proposed a class of double-loop approximation algorithms to iteratively approximate the stochastic gradient of the outer objective and incurred a sample complexity of  $\mathcal{O}(\eps^{-2})$ in order to achieve the $\eps$-stationary point. The double loop structure of these approaches made them impractical for large-scale problems; \cite{ghadimi2018approximation} required solving an inner optimization problem to a predefined accuracy, while \cite{yang2021provably} required a large batch size of $\mathcal{O}(\eps^{-1})$ at each iteration. To address this issue, various {single-loop methods, involving simultaneous update of inner and outer optimization variables, have been developed  \cite{hong2020two,chen2021single, khanduri2021near, yang2021provably}. A single-loop two-time scale stochastic algorithm proposed in \cite{hong2020two} incurred a sub-optimal sample complexity of $\mathcal{O}(\eps^{-2.5})$.} This is further improved recently in \cite{chen2021single, khanduri2021near, yang2021provably}, in which the authors have utilized the momentum-based variance reduction technique from \cite{cutkosky2019momentum} to obtain optimal convergence rates. While all of the above-mentioned works seek to solve \eqref{main_problem_2}, they are projection-based and require a projection on to $\cX$ at every iteration. In this work, we are interested in developing projection-free stochastic optimization algorithms for bi-level problems, which is still an open problem and the subject of the work in this paper. 

\noindent{\textbf{Compositional optimization}} problems have been recently studied and various algorithms have been proposed in  \cite{wang2017stochastic, wang2017accelerating,ghadimi2020single,chen2020solving,yuan2019stochastic,yang2020stochastic,thomdapu2020stochastic}. The  seminal work  in \cite{wang2017stochastic} proposed a quasi-gradient method called stochastic compositional gradient descent (CGD) to solve the problem via a two time-scale approach. In \cite{wang2017accelerating}, the authors proposed an accelerated SCGD method that achieved an improved sample complexity of $\mathcal{O}(\eps^{-2})$ for convex objectives and  $\mathcal{O}(\eps^{-2.25})$ for non-convex objectives. Further, different variance-reduced SCGD variants have been proposed, such as  SCVR \cite{liu2017variance}, VRSC-PG  \cite{huo2018accelerated}, SARAH-Compositional \cite{yuan2019efficient,yuan2019stochastic}, and STORM-Compositional \cite{yang2020stochastic}. In the literature, we can also find some single time-scale algorithms to solve compositional problems \cite{ghadimi2020single,chen2020solving}. Work in \cite{ghadimi2020single} presented a nested averaged stochastic approximation (NASA) and proved a sample complexity of  $\mathcal{O}(\eps^{-2})$. Recently,  another single time-scale algorithm called the stochastically corrected stochastic compositional gradient method (SCSC) is proposed in \cite{chen2020solving} that converges at the same rate as the SGD methods for non-compositional stochastic optimization. It further adopted the Adam-type adaptive gradient approach and achieved the optimal sample complexity of $\mathcal{O}(\eps^{-2})$. Again, all the above-mentioned algorithms either solve an unconstrained problem or use projection operation at each iteration to deal with the constraints. In this work, we developed a projection-free algorithm for compositional problems as well. Note that even-though compositional problems are a special case of bi-level problems, we have studied them separately in this work and proposed a novel projection-free algorithm specifically for compositional problems to achieve the optimal sample complexity. 

\noindent{\textbf{Projection-free algorithms}}  have been extensively studied to solve the single-level optimization  problems of the form \eqref{main_problem_2} in the literature  \cite{nemirovski2009robust,nitanda2014stochastic,lan2012optimal}. A number of first-order projection-free algorithms have been developed for stochastic optimization problems as well \cite{mokhtari2020stochastic,xie2020efficient,hazan2016variance,reddi2016stochastic,akhtar2021conservative}. The stochastic FW method proposed in \cite{hazan2016variance} achieves a sample complexity of $\mathcal{O}(\eps^{-3})$ but requires an the batch size  $b=\mathcal{O}(t)$, where $t$ is the iteration index. The need for increasing batch sizes was dropped in \cite{mokhtari2020stochastic}, which worked with a standard mini-batch but still achieved the same sample complexity. Finally, an improved stochastic recursive gradient estimator-based algorithm called ORGFW was proposed in \cite{xie2020efficient} and achieved a sample complexity of $\mathcal{O}(\eps^{-2})$. For non-convex problems, \cite{reddi2016stochastic} proposed an approach where the batch-size depends on the total number of iterations, resulting in a sample complexity of $\mathcal{O}(\eps^{-4})$. Later, work in \cite{zhang2020one} came up with a two-sample strategy and achieved  $\mathcal{O}(\eps^{-3})$ sample complexity. We remark that the idea of projection-free algorithms is limited only to single-level optimization problems in the existing literature. Therefore, there are no corresponding oracle complexity bounds for projection-free algorithms in bilevel settings, and our work fills this crucial gap.  We present our main contributions as follows.

\subsection{Contributions}
\begin{itemize}

\item Firstly, we propose a novel projection-free SBFW algorithm, utilizing the idea of momentum-based gradient tracking \cite{cutkosky2019momentum} in order to track the gradient of the outer objective function. The combination of this idea along with FW updates allows us to achieves the sample complexity of   $\mathcal{O}(\eps^{-3})$  and  $\mathcal{O}(\eps^{-4})$ for convex and non-convex cases, respectively (cf, Sec. \ref{sec:convergence}).  
 
\item 
 Secondly, in contrast to the existing literature,  we consider the compositional problems (which is a special case of bi-level problems) separately in this paper and propose a novel \textbf{S}tochastic \textbf{C}ompositional \textbf{F}rank \textbf{W}olfe (SCFW) algorithm. The SCFW algorithm is able to achieve a better convergence rate than bi-level problems under fewer assumptions on the inner objective functions. SCFW achieves the sample complexity of $\mathcal{O}(\eps^{-2})$ and  $\mathcal{O}(\eps^{-3})$ for convex and non-convex cases, respectively (cf, Sec. \ref{sec:convergence}). These results are interesting because they match with the existing best possible convergence
rates of projection-based methods (which are computationally expensive) for compositional problems \cite{wang2017accelerating,chen2020solving}, and projection-free methods for non-compositional (single-level) cases \cite{hassani2020stochastic,xie2020efficient}.
% Note that the proposed method is for projection-free compositional problems, which are different from the existing settings in the literature and is the first of their kind.	

\item
 Finally, we test the proposed algorithm on matrix completion and the problem of policy evaluation in reinforcement learning and establish the efficacy of the proposed techniques as compared to state-of-the-art algorithms (cf. Sec. \ref{sec:application}). We achieve an improvement of up to $82\%$ in the computation time for the proposed algorithm as compared to state-of-the-art methods.  
\end{itemize}
A comprehensive list of all existing related works is provided in Table \ref{table1}.{ For bi-level problems, the second column of Table \ref{table1} represents the objective type (convex/non-convex) of outer function. The inner objective is strongly convex for all the methods in  Table \ref{table1}.} 

\textbf{Notation:} First, we defined the compact notations we utilize in the convergence proofs. We denote column vectors with lowercase boldface $\x$, its transpose as $\x^\top$, and its Euclidean norm by $\norm{\x}$.  We use $\mbE_t:=\mbE[\cdot|\mathcal{F}_t]$ to denote
the conditional expectation  with respect to given sigma field $\mathcal{F}_t$
which contains all algorithm history (randomness) till step $t$.
\section{Problem Formulation}\label{sec:Pre}
Let us consider the optimization problem $\mathcal{P}_1$.  For most of the applications in practice \cite{rajeswaran2019meta,franceschi2018bilevel,zhang2020bi,hong2020two,borsos2020coresets}, the outer objective $F(\bbx, \y^\star(\x)):= \mbE_{\theta}[f(\y^\star(\x);\theta)] $ and the inner objective $G(\bby,\bbx) := \mbE_{\xi}[g(\y,\x;\xi)]$ are expected values of  continuous and proper closed functions $f:\mathbb{R}^m\rightarrow \mathbb{R}$, and the function $g$ is defined as $g(\cdot,\cdot,\xi):\mathbb{R}^n\times \mathbb{R}^m\rightarrow\mathbb{R}$ with respect to the {independent} random variables $\theta$ and $\xi$, respectively. Hence, the equivalent  bi-level stochastic optimization problem is given by
\begin{subequations}\label{first-level-sto}
	\begin{align}
	\x^*:=	&	\arg\min_{\x\in \cX \subset \mathbb{R}^m} Q(\x) = \mbE_{\theta}[f(\x,\y^\star(\x);\theta)]  ,
	\\
	&	\text{s.t.}\;\; \; \y^{\star}(\x)\in \arg \min_{\y \in \mathbb{R}^n}  \mbE_{\xi}[g(\y,\x;\xi)]. \label{second-level-sto}
	\end{align}
\end{subequations}
%
%Besides from the general form in \eqref{first-level-sto}, we are also interested in the special cases where it is possible to solve the inner objective for a given $\bbx$ and obtain a smooth closed-form expression for the optimal solution $\y^\star(\x)$.  
Besides the general form in \eqref{first-level-sto}, we are also interested in the special cases where it is possible to solve the inner optimization  problem at a given $\bbx$ and obtain a smooth closed form expression for the optimal solution $\y^\star(\x)$. For instance, if the inner objective $g$ is quadratic in $\y$, i.e. $g(\x,\y,\xi)=\norm{\y-h(\x,\xi)}^2$, where $h$ is a smooth function over $\x$, then we can write $\y^*(\bbx)=h(\x,\xi)$. Hence, the problem in \eqref{first-level-sto} boils down to a stochastic compositional  optimization problem given by
\begin{align} \label{compos}
\min_{\x\in \cX \subset \mathbb{R}^m} C(\x):=\mbE_{\theta} \left[f\left(\mbE_{\xi}[h(\x,\xi)],\theta \right)\right],&
\end{align}
which involves the nested expectations in the objective. The problem in \eqref{compos} has been independently considered in the literature and solved using two-time scale approaches \cite{wang2017stochastic,wang2016stochastic}. A single-time scale approach  for the problem in \eqref{compos} is also proposed in \cite{ghadimi2020single,chen2020solving}  but all the existing approaches are projection-based. In this work, we are interested in developing first-order methods to solve the problem in \eqref{first-level-sto} in a projection-free manner. We remark here that while the algorithms developed for problem in \eqref{first-level-sto} could be readily applied to solve the problems of the form in \eqref{compos}, there is a scope to further propose faster algorithms for compositional problems in \eqref{compos}. Therefore, we will consider the compositional problems separately from bi-level problems and derive better convergence rates. 
\begin{algorithm}[t]
	\caption{\textbf{S}tochastic \textbf{B}i-level \textbf{F}rank \textbf{W}olfe}
	\label{alg:SBFW}
%	\begin{algorithmic}
		\KwIn{$\x_1\in \cX, \y_1\in \mathbb{R}^m, \eta_t, \delta_t, \rho_t$, $\beta_t$, and $\d_1=h_1(\theta_1;\xi_1)\; ${using \eqref{biased_estimate}} }
		\For{$t=2$ {\bfseries to} $T$}
	{{\textbf{Update} approximate inner optimization solution 
		%\begin{align*}
		$$\y_{t}=\y_{t-1}-\delta_t \n_{\y} g(\x_{t-1},\y_{t-1},\xi_t)$$
		%\end{align*}
		\;
	 \textbf{Gradient tracking} evaluate $h(\x_t,\y_t;\theta_t,\xi_t)$ and $h(\x_{t-1},\y_{t-1};\theta_t,\xi_t)$ using \eqref{biased_estimate} and compute 
		\begin{align*}
		\d_t&=(1-\rho_t)(\d_{t-1}-h(\x_{t-1},\y_{t-1};\theta_t,\xi_t)\nonumber\\&\quad+ h(\x_t,\y_t;\theta_t,\xi_t)
		\end{align*}\;
	\textbf{Evaluate} feasible direction $\s_t=\arg \min_{\s\in \mathcal{X}}\ip{\s,\d_t}$\;
	 \textbf{Update} solution $\x_{t+1}=(1-\eta_{t})\x_t+\eta_{t} \s_t$}}
		\KwOut{$\x_{T+1}$ or $\hat{\x}$ selected uniformly from $\{\x_i\}_{i=1}^T$}
%	\end{algorithmic}
\end{algorithm}

\section{Algorithm Development}\label{Sec:Algo_Dev}
%\textbf{Solution and Challenges:}

In this section, we develop the proposed algorithm to solve the problem in \eqref{first-level-sto}.  We note that solving the bi-level optimization problem in \eqref{first-level-sto} is NP hard in general but we restrict our focus to problems where the  inner objective is continuously twice differentiable in $(\x, \y)$ and also strongly convex w.r.t $\y$ with parameter $\mu_g>0$. Such an assumption is common in the related works \cite{ghadimi2018approximation,chen2021single,yang2021provably} and ensures that for any $\x \in \cX$, $\y^\star(\x)$ is unique.  Applying any first order method to solve  \eqref{first-level-sto} requires the evaluation of the gradient of outer objective $\nabla Q(\x)$ with respect to $\bbx$. Let us denote the iteration index by $t\in\{1,2,\cdots,T\}$, and then we can write the standard projected gradient descent update as %
\begin{align}\label{projected_GD}
\bbx_{t+1}=\mathcal{P}_{\mathcal{X}}\left[\bbx_{t}-\alpha\nabla Q(\x_t)\right],
\end{align}
where $\mathcal{P}_{\mathcal{X}}[\cdot]$ denotes the projection onto the constraint set $\mathcal{X}$ and $\alpha$ is the constant step size. Note that the gradient calculation in \eqref{projected_GD} is achieved by application of \emph{implicit function theorem} to the optimality condition for inner optimization problem $\n_{\y}G(\bby,\bbx)=0$, calculating total derivative followed by chain rule to obtain the expression
\begin{align}\label{main_grad}
\nabla Q(\x_t)=&\n_\x F(\x_t,\y^\star(\x_t))-\n_{\x\y}^2G(\y^\star(\x_t), \x_t)\times[\n_{\y\y}^2G(\y^\star(\x_t), \x_t)]^{-1}\n_\y F(\x_t,\y^\star(\x_t)).
\end{align}
From \eqref{main_grad}, note that $\nabla Q(\x_t)$ requires the information about $\y^\star(\x_t)$ which is not available in general, unless the second level problem \eqref{second-level} has a closed form solution. Hence, it is not possible to utilize gradient based algorithms to solve the problem in \eqref{first-level-sto}. This challenge is addressed in the literature via the utilization of approximate gradients \cite{ghadimi2018approximation}. Following a similar approach, we use a surrogate gradient defined as $\n S(\x_t,\y_t)$ in place of original gradient $\n Q(\x_t)$. The surrogate gradient is obtained by replacing $\y^\star(\x_t)$ in \eqref{main_grad} with some $\y_t \in \mathbb{R}^n$ (we will define  the explicit value of $\y_t$ later) to write
\begin{align}\label{surrogate}
\n S(\x_t,\y_t)&=\n_\x F(\x_t,\y_t)-\n_{\x\y}^2G(\y_t, \x_t)\times[\n_{\y\y}^2G(\y_t, \x_t)]^{-1}\n_\y F(\x_t,\y_t).
\end{align}
\begin{algorithm}[t]
	\caption{\textbf{S}tochastic \textbf{C}ompositional \textbf{F}rank \textbf{W}olfe}
	\label{alg:SCFW}
%	\begin{algorithmic}
		\KwIn{$\x_0\in \cX, \eta_t, \delta_t, \rho_t$, $			\y_0=h(\x_0,\xi_0)$, $\d_0=\n F(\x_0,\y_{0},\thh_0,\xi_0)$}
		\For{$t=1$ {\bfseries to} $T$}
	 {\textbf{Update} inner function tracking \begin{align*}
		\y_{t}=(1-\delta_t)(\y_{t-1}-h(\x_{t-1},\xi_t))+ h(\x_t,\xi_t)
		\end{align*} \;
			 \textbf{Gradient tracking} \begin{align*}
		\d_t&=(1-\rho_t)(\d_{t-1}-\n F(\x_{t-1},\y_{t-1},\thh_t,\xi_t))+ \n F(\x_t,\y_{t},\thh_t,\xi_t)
		\end{align*}\;
	\textbf{Evaluate} feasible direction $\s_t=\arg \min_{\s\in \mathcal{X}}\ip{\s,\d_t}$\;
		\textbf{Update} solution $\x_{t+1}=(1-\eta_{t})\x_t+\eta_{t} \s_t$
}
		\KwOut{$\x_{T+1}$ or $\hat{\x}$ selected uniformly from $\{\x_i\}_{i=1}^T$}
%	\end{algorithmic}
\end{algorithm}
Even after replacing $\y^\star(\x_t)$ with some $\y_t \in \mathbb{R}^n$, there is an additional challenge of evaluating the individual terms in \eqref{surrogate} for the stochastic bi-level problems mentioned in \eqref{first-level-sto}. For instance, the term $\n_\x F(\x_t,\y_t)= \mbE_{\theta}[\n_\x f(\y_t;\theta)] $ in \eqref{surrogate} involves the evaluation of the expectation operator, which is not possible in practice due to the unknown data distribution.
One standard approach is to use unbiased stochastic gradient instead of the original gradient, but an unbiased estimate of $\n S(\x_t,\y\textbf{})$ would still require the computation of Hessian inverse. To avoid such complicated matrix computations, we follow the approach presented in \cite[Sec. E.4]{hong2020two} and compute a mini-batch approximation of Hessian inverse using the samples returned by the sampling oracle. {We assume the availability of sampling oracle such that for a given $\x \in \cX$ and $\y \in \mathbb{R}^n$, it returns unbiased samples $\n_\x f(\x,\y,\th)$, $\n_\y  f(\x,\y,\th)$, $\n_\y  g(\x,\y,\xi)$, $\n_{\x\y} ^2 g(\x,\y,\xi)$, and $\n_{\y\y} ^2 g(\x,\y,\xi)$  realized at random variables $\xi$ and $\thh$. 
	Having access to such oracle, we can write }the biased estimate (denoted by $h(\x_t,\y_t;\theta_t,\xi_t)$) of surrogate gradient $\n S(\x_t,\y_t)$ in \eqref{surrogate}  as 
\begin{align}\label{biased_estimate}
h(\x_t&,\y_t;\theta_t,\xi_t)=\n_\x f(\x_t,\y_{t};\theta_t)-M(\x_t,\y_t;\tilde{\xi}_t)\cdot\n_{\y}f(\x_t,\y_t;\theta_t),
\end{align}
where $\n_\x f(\x_t,\y_{t};\theta_t)$ is an unbiased estimate of $\n_\x F(\x_t,\y_t)$, $\n_\y f(\x_t,\y_{t};\theta_t)$ is an unbiased estimate of $\n_\y F(\x_t,\y_t)$, and  $M(\x_t,\y_t;\tilde{\xi}_t)$ is a biased estimate of product $\n_{\x\y}^2G(\y_t, \x_t)\cdot[\n_{\y\y}^2G(\y_t, \x_t)]^{-1}$. In the term $M(\x_t,\y_t;\tilde{\xi}_t)$, $\tilde{\xi}_{t}$ is defined as $\tilde{\xi}_{t}:=\{\xi_{t,i}: \: i\in\{0,1,\cdots,k\}\}$ which represents a collection of $(k+1)$ i.i.d. samples $\n^2_{\x\y}g(\y_t,\x_{t};\xi_{t,i})$. The explicit form of  $M(\x_t,\y_t;\tilde{\xi}_t)$ is
\begin{align}\label{M}
M(\x_{t},\y_{t};\tilde{\xi}_t)=&\n^2_{\x\y}g(\y_t,\x_{t};\xi_{t,0})\times\left[\frac{k}{L_g}\Pi_{i=1}^l\!\left(\!I\!-\!\frac{1}{L_g}\n_{\y\y}^2g(\y_t,\x_{t};\xi_{t,i})\right)\!\right],
\end{align}
where $l$ is selected uniformly from~$\{0,1,\cdots,k-1\}$.  Further, for $l$$=$$0$, we use the convention $$\Pi_{i=1}^l\left(I-\frac{1}{L_g}\n_{\y\y}^2g(\x_t,\y_{t};\xi_{t,i})\right)=I.$$ Hence, a stochastic version of the update in \eqref{projected_GD} would be given by
\begin{align}\label{projected_GD2}
\bbx_{t+1}=\mathcal{P}_{\mathcal{X}}\left[\bbx_{t}-\alpha  h(\x_t,\y_t;\theta_t,\xi_t)\right].
\end{align}
Similar to the update in \eqref{projected_GD2}, different variants are proposed in the literature \cite{ghadimi2018approximation,chen2021single,yang2021provably,khanduri2021near}. But a significant challenge that remains un-addressed to date in the literature for the bi-level problems is associated with the projection operator in \eqref{projected_GD2}. Projection-based algorithms require to perform a computationally expensive projection step at each iteration $t$. One projection step is called projection oracle call.  For instance, to achieve $\eps$ suboptimality, the projected subgradient
method \cite{beck2017first} entails $\mathcal{O}(\eps^{-2})$ projection oracle calls. The projection is easy to evaluate when the constraint set is a simple convex set (onto which projection operation is computationally cheap such as probability simplex) or has a closed-form solution (set of unit-ball). However, the projection step is often computationally costly (e.g., nuclear norm constraint), and its complexity could be comparable to the problem at hand \cite{jaggi2013revisiting}. 
In the next subsection, we obviate the issue related to projections by proposing projection-free algorithms for both the bi-level and compositional stochastic optimization problems, which is the key novel aspect of work in this paper.

% \textcolor{red}{The projection is easy to evaluate  when the constraint set is simple convex set or the dimension of $\bbx$ is small ($m$ is very small). But there are applications when the constraint set $\mathcal{X}$ is complicated to project on such as norm zero constraint (as in Example 1) or $m$ is really large (which is usually millions for federated learning applications as mentioned in Example 1). For instance, Can you provide an argument to motivate the complexity of the projection operator similar to the one given by Praneeth's paper.}

% The focus of this work paper is  on the case where the feasible set $\cX$ is complicated, in the sense that the computation cost of projection onto $\cX$ is high or even intractable. To be more precise, although the set $\cX$ is a closed convex set in $\mathbb{R}^n$ and is bounded with diameter $D$, i.e., $\norm{\x_1-\x_2} \leq D$ for all $\x_1,\x_2 \in \cX$, it does not allow for efficient projection  (for
% instance, projection on the nuclear norm ball). We further assume that the distribution of $\xi$ and $\thh$ are unknown, and we only have access to a stochastic oracle such that, for a given $\x \in \cX$, it returns a random vector $g(\x,\xi)$ and a noisy subgradient $\n g(\x,\xi)$ for some random sample $\xi$ and for a given vector $\y \in \mathbb{R}^m$, the oracle returns a noisy gradient $\n f(\y,\thh)$  for some random sample $\thh$.

%  
\subsection{Stochastic Projection-Free Bi-level Optimization Algorithm}
Before proceeding towards the design of a projection-free algorithm, a discussion regarding the particular choice of $\y_t$ in \eqref{projected_GD2} is due. A popular choice (see \cite{ghadimi2020single,chen2021single,hong2020two,ji2020bilevel}) for $\y_t$ is the stochastic gradient descent update for the inner optimization problem given by 
\begin{align}\label{inner_update}
\y_{t}=\y_{t-1}-\delta_t \n_{\y} g(\y_{t-1}, \x_{t-1};\xi_t), 
\end{align}
where $\n_{\y} g(\y_{t-1}, \x_{t-1};\xi_t)$ is the unbiased estimate of the gradient $\mathbb{E}_{\xi}\left[\n_{\y} g(\y_{t-1}, \x_{t-1};\xi)\right]$, and $\delta_t$ denotes the step size. 
%We remark that a similar updating rule for $\y_t$ as mentioned in \eqref{inner_update} is commonly used in the prior works, such as for inner optimization problems. 
Now we are ready to propose the first projection-free algorithm for the bi-level stochastic optimization problems. We propose to use a conditional gradient method (CGM) based approach instead of calculating the projection in \eqref{projected_GD2}. That is, we solve a  linear minimization problem to find a feasible direction $\bbs_t\in\mathcal{X}$ for a given stochastic gradient direction $ h(\x_t,\y_t;\theta_t,\xi_t)$, given by,
%
%\begin{align}
$\bbs_t:=\arg \min_{\s \in \cX} \ip{\s,  h(\x_t,\y_t;\theta_t,\xi_t)}$.
%\end{align}
%
This reduces the optimization problem of evaluating the projection operator in \eqref{projected_GD2} to solving a linear program which is easier to solve in practice. Hence, the iterate in \eqref{projected_GD2} gets modified to
\begin{align}
\bbs_t:=&\arg \min_{\x \in \cX} \ip{\x,  h(\x_t,\y_t;\theta_t,\xi_t)}\label{projected_GD22_1}
\\
\bbx_{t+1}=&(1-\eta_{t+1})\x_t+\eta_{t+1} \s_t,\label{projected_GD22_2}
\end{align}
where $\eta_{t}$ is the step size. To this end, we would like to emphasize that naive use of $h(\x_t,\y_t;\theta_t,\xi_t)$ in \eqref{projected_GD22_1} for the evaluation of $\bbs_t$ which is then used in \eqref{projected_GD22_2} can result in the iterate divergence due to the non-vanishing variance of the gradient estimate $h$. The standard approach to deal with this issue is to use a biased gradient estimate with low variance instead of an unbiased one. For example, a mini-batch approximation is proposed in  \cite{hazan2016variance,reddi2016stochastic} with linearly increasing batch size with iteration index. Such an approach runs into memory issues when utilized in practice. Another line of work suggests the use of gradient tracking given by $\d_t=(1-\rho_t)(\d_{t-1})+ \rho_t h(\x_t,\y_t;\theta_t,\xi_t)$, where $\rho_t$ being the tracking parameter, which does not suffer from the problem of increasing batch size. But such gradient tracking schemes are shown to be suboptimal even for single-level optimization problems \cite{mokhtari2020stochastic}, so no scope for much harder bi-level problems are considered in this work. To address the issue of memory and  iterate divergence, we took motivation from the momentum-based approach in \cite{cutkosky2019momentum} and propose to use the following gradient tracking scheme given by
\begin{align}\label{momentum-track}
\d_t=&(1-\rho_t)(\d_{t-1}-h(\x_{t-1},\y_{t-1};\theta_t,\xi_t))+ h(\x_t,\y_t;\theta_t,\xi_t).
\end{align}
We remark that such a tracking technique is recently utilized in \cite{khanduri2021near} for projection-based bi-level optimization problems. However, in this work, our focus lies in developing projection-free algorithms, and hence analysis is significantly different from \cite{khanduri2021near}. We will show in Lemma \ref{Lemma:d:SBFW} ( supplementary material) that momentum-based tracking technique such as in \eqref{momentum-track} results in a reduced variance for the gradient estimate, hence eventually resulting in improving the convergence rate of the algorithm to achieve the optimality guarantees. Now we proceed to modify the updates in \eqref{projected_GD22_1}-\eqref{projected_GD22_2} in order to propose the main projection-free updates as  
\begin{align}
\bbs_t:=&\arg \min_{\x \in \cX} \ip{\x,  \d_t}\label{projected_GD22_3}
\\
\bbx_{t+1}=&(1-\eta_{t+1})\x_t+\eta_{t+1} \s_t.\label{projected_GD22_4}
\end{align}
The proposed algorithm is summarized in Algorithm \ref{alg:SBFW}. The
details about the specific choice of step sizes $\delta_t$, $\rho_t$, and $\eta_{t}$ are discussed in Sec. \ref{sec:convergence}.

%
%
%%%%%%%%
%%%%%%%
%%START PROPOSED SCFW
%%%%%%%%%
%%%%%%%%%%%%
\subsection{Stochastic Projection-Free Compositional Optimization Algorithm}  \label{Sec:SCGD}
In this section, we propose a separate novel algorithm for problems of the form in \eqref{compos} in a projection-free manner. To solve compositional problems in \eqref{compos}, the classical SGD methods requires the unbiased samples of the actual gradient $\n C(\x_t)$  given by $\n C(\x_t,\thh_t,\xi_t):=\n h(\x_t, \xi_t)^\top\n f(\mbE_{\xi}[h(\x_t,\xi)],\thh_t)$. That is, for a given vector $\x_t$ and random sample $\xi_t$ and $\thh_t$, this gradient computation requires $\mbE_{\xi}[h(\x_t,\xi)]$ which is not available on a single query to sampling oracle. The existing projection-based stochastic compositional methods adopt some kind of tracking approach to approximate $\mbE_{\xi}[h(\x_t,\xi)]$ using the samples returned by the oracle. For instance, the seminal work SCGD in \cite{wang2017stochastic} approximated $h(\x_t)$ at each iteration as 
$\y_{t}=(1-\delta_t)(\y_{t-1})+ \delta_t h(\x_t,\xi_t)
$ 
where $\delta_t$ is a diminishing step size but achieved suboptimal convergence rate of $\mathcal{O}(1/\epsilon^{4})$. Recently, work in \cite{chen2020solving} proposed Adam-SCSC that  introduces a stochastic correction to the original SCGD \cite{wang2017stochastic} using Adam-type adaptive gradient approach and achieved the convergence rate of $\mathcal{O}(1/\epsilon^{2})$. However, all these existing  algorithms require projection onto the set $\cX$ at each iteration. In the proposed algorithm SCFW, slightly different from SBFW, we proposed to use momentum-based tracking technique to estimate $\bby_t$ as 
\begin{align}\label{inner_fucn_track}
\y_{t}=(1-\delta_t)(\y_{t-1}-h(\x_{t-1},\xi_t))+ h(\x_t,\xi_t).
\end{align} 
After utilizing $\y_t$, we denote the stochastic gradient estimate of compositional objective as $\n C(\x_t,\y_{t},\thh_t,\xi_t)=\n h(\x_t, \xi_t)^\top\n f(\y_t,\thh_t)$, and then similar to SBFW, we proposed to use the momentum based method to track the compositional gradient as  
%gradient $\n F(\x,\y)=\mbE[\n F(\x,\y,\thh,\xi)]$. That is, for the given i.i.d. sample $\{\xi_{t},\theta_t\}$,
% we estimate the required gradient using a
%stochastic recursive estimator as
%
\begin{align}\label{grad_track}
\d_t&=(1-\rho_t)(\d_{t-1}-\n C(\x_{t-1},\y_{t-1},\thh_t,\xi_t))+ \n C(\x_t,\y_{t},\thh_t,\xi_t).
\end{align}
We will establish that using such momentum based tracking technique for both the inner function \eqref{inner_fucn_track} as well as on the gradient of the objective \eqref{grad_track} not only reduces the variance of the approximation noises with each iteration but also helps in establishing optimal convergence results. The rest of the steps essentially remains the same as SBFW; that is, we solve a linear minimization problem $\s_t=\arg \min_{\s\in \mathcal{X}}\ip{\s,\d_t}$ and then update the iterate as $\x_{t+1}=(1-\eta_{t+1})\x_t+\eta_{t+1} \s_t$.
The proposed SCFW algorithm is summarized in Algorithm \ref{alg:SCFW}. The
details about the specific choice of step sizes $\delta_t$, $\rho_t$ and $\eta_{t}$ are  provided
in Sec. \ref{sec:convergence}.

%
%%%%%%%%
%%%%%%%
%%END PROPOSED SCFW
%%%%%%%%%
%%%%%%%%%%%%
%
%
%
%
%

\section{Convergence Analysis}\label{sec:convergence}
This section presents the convergence rate analysis for the proposed Algorithm \ref{alg:SBFW} and Algorithm \ref{alg:SCFW}. We provide the convergence rate results for both convex and non-convex objectives. Note that the inner objective is considered to be strongly convex for both cases.   Before proceeding towards the analysis, we first discuss the convergence criteria we utilize to evaluate the performance of the proposed algorithms.

\textbf{Convergence Criteria:} For the convex objective $Q(\bbx)$, we use the expected suboptimality $\mbE[Q(\x_T)-Q(\x^\star)]$ after $T$ number of iterations where $\x^*$ is as defined in \eqref{first-level-sto}. The expectation here is with respect to the randomness in the objective as well as iterate $\bbx_T$. However, for non-convex $Q(\x)$, we use the expected Frank-Wolfe gap $\mbE[\mathcal{G}(\x)]$ as the performance metric where $\mathcal{G}(\x)$ is defined as
\begin{align}
\label{conv-criteria:non convex}
\mathcal{G}(\x):=\max_{\mathbf{v}  \in \mathcal{X}} \ip{\mathbf{v}-\x,-\n Q(\x)}. 
\end{align}
The Frank-Wolfe gap is a standard performance metric for the constrained non-convex settings as mentioned in \cite{zhang2020one,reddi2016stochastic,lacoste2016convergence}. From the definition in \eqref{conv-criteria:non convex}, we note that $\mathcal{G}(\x)\geq 0$ for all $\x\in \cX$ and if $\exists\x'\in \cX$ such that $\mathcal{G}(\x')= 0$, then $\x'$ is a first-order stationary point. { Further, for comparison with other methods, we will use the stochastic first-order oracle (SFO) complexity which is a commonly used metric to compare stochastic first-order methods. SFO  defines the total number of times an algorithm is required to call a first-order oracle (which provides stochastic gradients) to reach $\eps$-approximate (or stationary) solution. Another metric that is often used to evaluate stochastic projection-free algorithms is the linear minimization oracle (LMO) complexity, which is the total number of times an algorithm needs to solve the linear minimization problem to reach $\eps$-approximate (or stationary) solution. However, since both the proposed SBFW and SCFW algorithms require $\mathcal{O}(1)$ LMO calls at each iteration, the LMO complexity is similar to the SFO complexity. We will discuss the SFO complexity of each algorithm after the main theorems.}

Next, we state the assumptions required to perform the analysis in this work that is similar to the assumptions considered in the existing literature \cite{hong2020two,khanduri2021momentum}.   
\begin{assumption}\label{ass:SBFW} For the bi-level stochastic optimization problem in \eqref{first-level-sto}, we need  the following statements to hold true for the analysis. 
		\begin{enumerate}[(i)]
	%		\item \label{assSBFW:sampling_oracle}	\normalfont (\emph{Sampling oracle})
	%		For a given $\x \in \cX$ and $\y \in \mathbb{R}^n$, the oracle return samples $\n_\x f(\x,\y,\th)$, $\n_\y  f(\x,\y,\th)$, $\n_\y  g(\x,\y,\xi)$, $\n_{\x\y} ^2 g(\x,\y,\xi)$ and $\n_{\y\y} ^2 g(\x,\y,\xi)$  realized at random variables $\xi$ and $\thh$. All of these samples are unbiased which implies that $\mbE[\n_\x f(\x,\y,\th)]=\n_\x f(\x,\y)$, $\mbE[\n_\y  f(\x,\y,\th)]=\n_\y  f(\x,\y)$, $\mbE[\n_\y  g(\x,\y,\xi)]=\n_\y  g(\x,\y)$, $\mbE[\n_{\x\y} ^2 g(\x,\y,\xi)]=\n_{\x\y} ^2 g(\x,\y)$, and $\mbE[\n_{\y\y} ^2 g(\x,\y,\xi)]=\n_{\y\y} ^2 g(\x,\y)$.
	
\item \label{assSBFW:bounded variance}
	The stochastic gradient estimates satisfies 
	\begin{align}
	&	\mbE[\norm{\n_\x f(\x,\y)-\n_\x f(\x,\y,\theta)}^2]\leq \sigma^2_{\x}, \nonumber\\
	&	\mbE[\norm{\n_\y  f(\x,\y)-\n_\y  f(\x,\y,\theta)}^2]\leq \sigma^2_{\y},\nonumber\\
	&	\mbE[\norm{\n_{\x\y}^2 g(\x,\y)-\n_{\x\y}^2  g(\x,\y,\theta)}^2]\leq \sigma^2_{\x\y},\nonumber\\
	&	\mbE[\norm{\n_\y g(\x,\y)-\n_\y g(\x,\y,\xi)}^2]\leq \sigma^2_{g}\nonumber,
	\end{align}
	for some  $\sigma^2_{\x}>0$, $\sigma^2_{\y}>0$, $\sigma^2_{\x\y}>0$, and $\sigma^2_{g}>0$. 
	
\item \label{assSBFW:Lipschitz_x} For any given $\x \in \cX$,  the terms $\n_\x f(\x,\y)$, $\n_\y f(\x,\y)$, $\n_\y  g(\x,\y)$, $\n_{\x\y} ^2g(\x,\y)$ and $\n_{\y\y} ^2g(\x,\y)$ are Lipschitz continuous with respect to $\y$ with Lipschitz parameter  $L_{f_\x}$, $L_{f_\y}$, $L_{g}$, $L_{g_{\x\y}}$ and $L_{g_{\y\y}}$, respectively. 
	
\item \label{assSBFW:Lipschitz_y} For any given $\y \in \mathbb{R}^n$, the terms $\n_\y  f(\x,\y)$,  $\n_{\x\y} ^2g(\x,\y)$ and $\n_{\y\y} ^2g(\x,\y)$ are Lipschitz continuous with respect to $\x$ with positive constants  $L_{f_\y}$, $L_{g_{\x\y}}$ and $L_{g_{\y\y}}$, respectively. Note that for the sake of simplicity, here we slightly abused the notation and used the same constants as in Assumption \ref{ass:SBFW} (ii).

\item \label{assSBFW:bounded_moment}
	For all $\x\in \cX$ and $\y \in \mathbb{R}^n$, it holds that $\mbE[\norm{\n_\y  f(\x,\y)}]\leq C_{\y}$ and 
	$\mbE[\norm{\n_{\x\y} ^2 g(\x,\y)}]\leq C_{\x\y}$ for some for constants $C_{\y}>0$ and   $C_{\x\y}>0$.

	(v) \label{assSBFW:strong_convexity}
	The inner function $g(\x,\y)$ is $\mu_g$-strongly convex in $\y$ for any $\x\in \cX$.
	
		\end{enumerate}
\end{assumption}

For SCFW, the basic assumptions are essentially the same as SBFW; however, for the sake of clarity, we state them here explicitly. It must be noted that for compositional problems, the inner function $h(\cdot)$ is not necessarily required to be strongly convex. 
\begin{assumption}\label{ass:SCFW}
	For the compositional stochastic optimization problem in \eqref{first-level-sto}, we need the following statements to hold true for the analysis.

	\begin{enumerate}[(i)]
		\item \label{ass:sampling_oracle} (\emph{Sampling oracle})
		For a given $\x \in \cX$ and $\y \in \mathbb{R}^n$, the oracle returns samples $h(\x,\xi)$, $\n h(\x,\xi)$ and $\n f(\y,\thh)$  for some random sample $\xi$ and $\thh$. These samples are unbiased that is  $\mbE_\xi[h(\x,\xi)]=h(\x)$ and $\mbE_{\{\xi,\thh\}}[\n h(\x,\xi)^\top \n f(\y,\thh)]=\n h(\x)^\top \n f(\y)$.

		\item \label{ass:inner_var}
		 (\emph{Bounded Variance})
		The inner  function $h(\cdot)$ has bounded variance, i.e., 
		\begin{align}
		\mbE[\norm{h(\x,\xi)-h(\x)}]\leq \sigma_h^2.
		\end{align}
		
		\item \label{ass:smoothness}
		 (\emph{Lipschitz continuous})
		The functions $f(\cdot)$ and $h(\cdot)$ are smooth, i.e., for any $\y_1,\y_2 \in \mathbb{R}^n$
		\begin{align}
		\norm{\n f(\y_1,\thh)-\n f(\y_2,\thh)}\leq L_f \norm{\y_1-\y_2},
		\end{align}
		and for any $\x_1,\x_2 \in \cX$, it holds that
		\begin{align}
		\norm{\n h(\x_1,\xi)-\n h(\x_2,\xi)}\leq L_h \norm{\x_1-\x_2}.
		\end{align}
		
		\item \label{ass:bounded_moment}
		 (\emph{Bounded second moments})
		The stochastic gradients of functions $f(\cdot)$ and $h(\cdot)$ have bounded second order moments, i.e., for any $\x\in \cX$, and $\y \in \mathbb{R}^n$
		\begin{align}
		\mbE[\norm{\n f(\y,\thh)}^2]\leq M_f;\hspace{0.3cm}
		\mbE[\norm{\n h(\x,\xi)}^2]\leq M_h.
		\end{align}

	\end{enumerate}
	
\end{assumption}
Note that for SCFW, we have assumed both the inner and outer functions to be smooth; the resulting composition function $C(\x)$ will also be smooth \cite{zhang2019stochastic}. Its  smoothness parameter $L_F$ can be easily obtained using Assumption \ref{ass:SCFW} (\ref{ass:smoothness}-\ref{ass:bounded_moment}) as $L_F=M_h^2L_f + M_fL_h$.

We start the analysis by presenting intermediate Lemmas \ref{Lemma: ySBFW}-\ref{Lemma:d:SBFW} and Corollary \ref{coro1} which eventually leads to the main result of this section presented in Theorem \ref{TH:convex-bilevel}. 
\begin{lemma}\label{Lemma: ySBFW}  Consider the proposed Algorithm \ref{alg:SBFW} and  $\x_t$ be the iterates generated by it, then for the algorithm parameter $\delta_t\leq \min\{\frac{2}{3\mu_g},\frac{\mu_g}{2(1+\sigma_g^2)L_g^2}\}$ and step size $\eta_t$, the optimality gap of the lower level problem satisfies 
	\begin{align}
	\mbE_t[\norm{\y_t-\y^\star(\x_t)}^2]&\leq \left(1-\frac{\delta_t\mu_g}{2}\right)\mbE_t[\norm{\y_{t-1}-\y^\star(\x_{t-1})}^2]\nonumber\\&\quad+\frac{2\eta_{t-1}^2}{\delta_t\mu_g}\left(\frac{C_{\x\y}}{\mu_g}\right)^2D^2+4\delta_t^2\sigma_g^2.
	\end{align}
\end{lemma}
The proof Lemma \ref{Lemma: ySBFW} is provided in section \ref{proof:Lemma: ySBFW} of the supplementary material. Lemma \ref{Lemma: ySBFW} quantifies how close $\y_t$
is from the optimal solution of inner problem at $\x_t$ and establishes the progress of the inner-level update.
%\begin{lemma}\label{old_lemma}[Lemma 11, \cite{hong2020two} ]
%	Consider the biased estimate defined in \eqref{biased_estimate}. Then for any $k\geq 1$, under Assumption \ref{ass:SBFW}, it hold that 
%	\begin{align}
%	&\!\!\!\!\norm{\n S(\x,\y)\!-\!\mbE[h(\x,\y;\theta,\xi)]}\!\leq\!\! \frac{C_{\x\y}C_{\y}}{\mu_g}\left(1\!-\!\frac{\mu_g}{L_g}\right)^k\!\label{old1}
%	\\
%	&\!\!\!\! \mbE[\|h(\x_t,\y_t;\theta_t,\xi_t)-\mbE[h(\x_t,\y_t;\theta_t,\xi_t)]\|^2]\leq \sigma_f^2,\label{old2}
%	\end{align}
%	where, $\sigma_f^2$$=$$\sigma_{\x}^2+\frac{3}{\mu_g^2}\left[(\sigma_{\y}^2+C_{\y}^2)(\sigma_{\x\y}^2+2C_{\x\y}^2)+\sigma_{\y}^2C_{\x\y}^2\right]$.
%\end{lemma}
%the specific choice of    mentioned in  satisfies the Assumption . In fact for $\beta$ of the order $\mathcal{O}((1-\frac{\mu_g}{L_g})^k)$, it can be shown that one can construct $h(\x_t,\y_t;\theta_t,\xi_t)$ with just $\mathcal{O}(1+log\; t)$ calls to the stochastic oracles while satisfying the requirement in Assumption \eqref{assSBFW:biased estimate} [\cite{hong2020two} Lemma 11]. Here, $k$ is the number of samples used for
%
%

\begin{lemma}\label{Lemma:d:SBFW}  Consider the proposed Algorithm \ref{alg:SBFW} and  $\x_t$ be the iterates generated by it, then for the algorithm parameter $\delta_t$, $\rho_t$ and $\eta_t$, we have 
	\begin{align}\label{ooo1}
	&\mbE_t[\|\d_t-\n S(\x_t,\y_{t})-B_t\|^2\nonumber\\&\leq (1-\rho_t)^2\mbE_t[\|(\d_{t-1}-\n S(\x_{t-1},\y_{t-1})-B_{t-1})\|^2]\nonumber\\&\quad+4L_k\delta_t^2L_g^2+4L_k\eta_{t-1}^2D^2+2\rho_t^2\sigma_f^2.
	\end{align}
\end{lemma}
The proof of Lemma \ref{Lemma:d:SBFW} is provided in section  \ref{Proof:Lemma:dSBFW} of the supplementary material.

Lemma \ref{Lemma:d:SBFW} describe the tracking error in the gradient approximation $\n S(\x,\y)$ at point $\x_t$ and $\y_t$. The presence of $(1-\rho_t)^2$ term in RHS of \eqref{ooo1} shows that the variance of the tracking error reduces with iteration.

Next,  we utilize Lemma \ref{Lemma: ySBFW}-\eqref{Lemma:d:SBFW} to establish a bound on the gradient estimation error for $\n Q$ presented in the form of following Corollary \ref{coroSBFW}.

\begin{corollary}\label{coroSBFW} For the proposed Algorithm \ref{alg:SBFW}, with $\delta_t=\frac{2a_0}{t^{q}}$,  where $a_0=\min\{\frac{1}{3\mu_g},\frac{\mu_g}{2(1+\sigma_g^2)L_g^2}\}$, $\rho_t=\frac{2}{t^{q}}$, $\beta_t\leq \frac{C_{\x\y}C_{\y}}{\mu_g(t+1)^q}$ and $\eta_t\leq\frac{2}{(t+1)^{3q/2}}$ for $0<q\leq 1$,  the gradient approximation error $\mbE\norm{\n Q(\x_t)-\d_t}^2$ converges to zero at the following rate
	\begin{align}
	\mbE\norm{\n Q(\x_t)-\d_t}^2 \leq \frac{C_1}{(t+1)^{q}},
	\end{align}
	where $C_1=3(\max\{2^q\norm{\y_1-\y^\star(\x_1)}^2,(2(C_{\x\y}/\mu_g)^2D^2$ $+16a_0^2\sigma_g^2)/(2a_0-1)\}+\tfrac{C_{\x\y}C_{\y}}{\mu_g}+8(2L_kL_g^2+L_kD^2+\sigma_f^2))$.
\end{corollary}
The proof of Corollary \ref{coroSBFW} is provided in section  \ref{Proof:corollary:SBFW} of the supplementary material.
The result in Corollary \ref{coroSBFW} is presented in general form and indicates that for properly chosen parameters $q$, the gradient
approximation error in expectation decreases at each iteration and approaches zero asymptotically. We will use this upper bound to prove the convergence of the proposed
algorithm SBFW for a different types of objective functions in the following theorem. Note that in the analysis of Corollary \ref{coroSBFW} we have set $\beta_t\leq \frac{C_{\x\y}C_{\y}}{\mu_g(t+1)^q}$. To satisfy this condition, the number of samples $k$ at iteration $t$ needed to
approximate the Hessian inverse in \eqref{M} is  $k=\mathcal{O}(\log((1+t)^q))$.

Now we are ready to present the first main result of this work as Theorem \ref{TH:convex-bilevel}.
%\subsection{Convergence Rate of SBFW (Algorithm \ref{alg:SBFW})}\label{conv:SBFW}
%
\begin{thm}\label{TH:convex-bilevel}
	(\textbf{Convergence Rate of SBFW}) Consider the proposed Algorithm \ref{alg:SBFW} and suppose Assumption \ref{ass:SBFW} is satisfied. Then,
	
	\noindent(i) (Convex Bi-Level): If $Q$ is convex on $\cX$ and we set  $\delta_t=\frac{a_0}{t^{2/3}}$,  where $a_0=\min\{\frac{2}{3\mu_g},\frac{\mu_g}{2(1+\sigma_g^2)L_g^2}\}$, $\rho_t=\frac{2}{t^{2/3}}$, $\eta_t=\frac{2}{t+1}$ and $k=\frac{2L_g}{3\mu_g}(\log(1+t))$, then the output is feasible $\x_{T+1}\in \cX $ and satisfies
	\begin{align}
	\mbE[Q(\x_{T+1})-Q(\x^\star)]
	&\leq \frac{12 D\sqrt{C_1}}{5(T+1)^{\frac{1}{3}}}+\frac{2L_QD^2}{(T+1)}.
	\end{align}
	
	\noindent (ii) (Non-Convex Bi-level): If $Q$ is non-convex and we set  $\delta_t=\frac{a_0}{t^{1/2}}$,  where $a_0=\min\{\frac{2}{3\mu_g},\frac{\mu_g}{2(1+\sigma_g^2)L_g^2}\}$, $\rho_t=\frac{2}{t^{1/2}}$, $\eta_t=\frac{2}{(T+1)^{3/4}}$ and $k=\frac{L_g}{2\mu_g}(\log(1+t))$ , then the output is feasible $\hat{\x} \in \cX $ and satisfies
	\begin{align*}
	&\mbE[\mathcal{G}(\hat{\x})]\leq \frac{Q(\x_1)-Q(\x^\star)}{(T+1)^{1/4}}+\frac{16 D\sqrt{C_1}}{3(T+1)^{1/4}}+\frac{L_QD^2}{(T+1)^{3/4}},
	\end{align*}		
	here $L_Q=\frac{(L_{f_y}+L)C_{\x\y}}{\mu_g}+L_{f_x}+C_y\left[\frac{L_{g_{\x\y}}C_y}{\mu_g}+\frac{L_{g_{yy}}C_xy}{\mu_g^2}\right]$ and $C_1=3(\max\{2\norm{\y_1-\y^\star(\x_1)}^2,(2(C_{\x\y}/\mu_g)^2D^2+16a_0^2\sigma_g^2)/(2a_0-1)\}+\tfrac{C_{\x\y}C_{\y}}{\mu_g}+8(2L_kL_g^2+L_kD^2+\sigma_f^2))$.
\end{thm}
The proof of Theorem \ref{TH:convex-bilevel} is provided in Appendix \ref{proof:theorem_convex_bilevel} in the supplementary material. Theorem \ref{TH:convex-bilevel} shows that the optimality gap for SBFW decays as $\mathcal{O}(T^{-1/3})$ for general convex objectives and for non-convex case it establishes an upper bound on the expected Frank-Wolfe gap for the iterates generated by SBFW that converges to zero
at least at the rate of $\mathcal{O}(T^{-1/4})$, where $T$ is the total number of iterations. It must be noted that for at each iteration, SBFW requires $2k+1$ gradient samples to obtain gradient estimate: $2k$ samples for outer gradient estimate \eqref{momentum-track} and one sample for inner variable update \eqref{inner_update}. Further, we have set $k\approx\mathcal{O}(\log(t))$. Hence, the SFO complexity of SBFW for outer objective is $\mathcal{O}(\log(\eps^{-1})\eps^{-3})\approx \mathcal{O}(\eps^{-3}) $ and $\mathcal{O}(\log(\eps^{-1})\eps^{-4})\approx \mathcal{O}(\eps^{-4}) $ for convex and non-convex objective, respectively. Similarly, observe that $\mbE[\norm{\y_{t}-\y^\star(\x_t)}]\leq \mathcal{O}( (t+1)^{-q})$ (see \eqref{final} in the supplementary material), where $q=2/3$ (as $\delta_t=\mathcal{O}(t^{-2/3})$) for convex objective and $q=1/2$ (as $\delta_t=\mathcal{O}(t^{-1/2})$) for non-convex objective. Hence,  the SFO complexity of inner objective for SBFW is $\mathcal{O}(\eps^{-1.5})$ and $\mathcal{O}(\eps^{-2})$ for convex and non-convex  function, respectively. It can be seen that complexity for inner level objective of the proposed algorithm SBFW is comparable to the projection-based state-of-the-art methods \cite{ghadimi2018approximation,hong2020two,ji2020bilevel,chen2021single,khanduri2021near}, however, it shows slightly worse performance in terms of the outer level complexity. This is not surprising as we are tackling the outer level in a projection-free manner. 

Next in order to study the convergence rate of the proposed algorithm SCFW, we state some results on the error of the inner function approximation and gradient approximation in the form of following lemma.
\begin{lemma}\label{Lemma:ySCFW}  Consider the proposed Algorithm \ref{alg:SCFW} and  $\x_t$ be the iterates generated by it, then the sequence $\norm{\y_{t}-h(\x_{t})}^2$ converges to zero at the following rate
	\begin{align}\label{Lemma: y eq}
	\mbE\norm{\y_{t}-h(\x_t)}^2&\leq    (1-\delta_t)^2\mbE\norm{(\y_{t-1}-h(\x_{t-1}))}^2+2\delta_t^2\sigma_h^2\nonumber\\&\quad+2(1-\delta_t)^2\eta_{t-1}^2M_hD^2.
	\end{align}	
\end{lemma}

%The proof of Lemma \ref{Lemma:ySCFW}  

The proof of Lemma \ref{Lemma:ySCFW} is provided in Appendix \ref{proof:Lemma:ySCFW}. It shows that the distance between $\y_t$ and $h(\x_t)$ decreases with iteration in expectation.  Intuitively, this means that our tracking variable $\y_t$ will converge to the unknown $h(\x_t)$. This result will be used to obtain the convergence rates of the proposed SCFW algorithm for different kinds of objective functions, which we provide as Theorem \ref{theorem1}.

Note that as instead of directly using the gradient samples, we are approximating it using tracking technique, we provide an upper bound on the gradient approximation error as follows.
\begin{corollary}\label{coro1} For the proposed Algorithm \ref{alg:SCFW}, with $\delta_t=\frac{2}{t^{p}}$, $\rho_t=\frac{2}{t^{p}}$ and $\eta_t\leq\frac{2}{(t+1)^{p}}$ for $0<p\leq 1$,  the gradient approximation error $\mbE\norm{\n C(\x_t)-\d_t}^2$ converges to zero at the following rate
	\begin{align}
	\mbE\norm{\n C(\x_t)-\d_t}^2 \leq \frac{A_1}{(t+1)^{p}}
	\end{align}
	where $A_1:=32[M_h(M_f+28L_f\sigma_h^2)+(M_f L_h+28M_h^2L_f)D^2]$.
\end{corollary}
The proof of Corollary \ref{coro1} in section \ref{proof:coro1} of the supplementary material. The result in Corollary \ref{coro1} indicates that for properly chosen parameters $p$, the gradient
approximation error in expectation decreases at each iteration and approaches zero asymptotically. We will use this upper bound to prove the convergence of the proposed
algorithm SCFW, which we will discuss next.

% Here, $k$ is the number samples required for approximation defined in \eqref{M}.
%
%%%%%%%%%%%%%%%%%%%%%%
%%%%%%%%%%%%%%%%%%%%%%%%%
%%%%%%%%%%%%%%%%%%%%%%%%%%%
% START: CONVERGENCE ANA OF SCFW
%%%%%%%%%%%%%%%%%%%%%%%%%%%%
%%%%%%%%%%%%%%%%%%%%%%%%%%
%\subsection{Convergence Analysis of SCFW (Algorithm \ref{alg:SCFW})}\label{conv:SCFW}

\begin{thm}\label{theorem1}
	(\textbf{Convergence Rate of SCFW}) Consider the proposed Algorithm \ref{alg:SCFW} and suppose that  Assumption \ref{ass:SCFW} is satisfied. Then, 
	
	\noindent (i)
	(Convex Compositional): If $C$ is convex on $\cX$ and we set $\rho_t=\frac{2}{t}$, $\delta_t=\frac{2}{t}$ and $\eta_t=\frac{2}{t+1}$, then the output is feasible $\x_{T+1}\in \cX $ and satisfies
	\begin{align}
	\mbE[C(\x_{T+1})-C(\x^\star)]
	&\leq  \frac{8 D\sqrt{A_1}}{3(T+1)^{\frac{1}{2}}}+\frac{2L_FD^2}{(T+1)},
	\end{align}
	
	\noindent (ii) (Non-Convex Compositional): If $C$ is non-convex and we set $\rho_t=\frac{2}{t^{2/3}}$, $\delta_t=\frac{2}{t^{2/3}}$ and $\eta_t=\frac{2}{(T+1)^{2/3}}$, then the output is feasible $\hat{\x} \in \cX $ and satisfies
	\begin{align*}
	\mbE[\mathcal{G}(\hat{\x})]\leq \frac{C(\x_1)-C(\x^\star)}{(T+1)^{1/3}}+\frac{6 D\sqrt{A_1}}{(T+1)^{1/3}}+\frac{L_FD^2}{(T+1)^{2/3}},
	\end{align*}
	where $L_F=M_h^2L_f + M_fL_h$ and $A_1:=32[M_h(M_f+28L_f\sigma_h^2)+(M_f L_h+28M_h^2L_f)D^2]$.		
\end{thm}
The proof of Theorem \ref{theorem1} is provided in section \ref{proof:theorem1} of the supplementary material. Interestingly, the optimality gap decays as $\mathcal{O}(T^{-1/2})$, which is
the optimal rate even for projected stochastic compositional optimization problems with general convex objectives \cite{wang2017accelerating}. 
For the non-convex case, it establishes an upper bound on the expected Frank-Wolfe gap for the iterates generated by SCFW and shows that it converges to zero
at least at the rate of $\mathcal{O}(T^{-1/3})$, where $T$ is the total number of iterations. It must be noted that at each iteration, SCFW requires only two gradient samples to obtain gradient estimate. Hence, it has SFO complexity is $\mathcal{O}(\eps^{-2})$ and $\mathcal{O}(\eps^{-3})$ for convex and non-convex objectives, respectively. Interestingly, these results match with the state-of-the-art methods \cite{zhang2020one,xie2020efficient} for projection-free single level (or non-compositional) stochastic  optimization.

\textit{Remark:} Note that  SCFW analysis is not a straightforward extension to SBFW, this is because we specifically re-derived inequality to exploit the special case of SCFW where we have an additional tracking step for inner function ($h(\textbf{x},\xi)$ of eq. \eqref{compos}), hence requiring a separate analysis.  In fact, note that we are able to achieve better convergence rate for SCFW $(\mathcal{O}(\epsilon^{-2}),\mathcal{O}(\epsilon^{-3}))$ as compared to SBFW $(\mathcal{O}(\epsilon^{-3}),\mathcal{O}(\epsilon^{-4}))$. 

%%%%%%%%%%%%%%%%%%%%%%
%%%%%%%%%%%%%%%%%%%%%%%%%
%%%%%%%%%%%%%%%%%%%%%%%%%%%
% END: CONVERGENCE ANA OF SCFW
%%%%%%%%%%%%%%%%%%%%%%%%%%%%
%%%%%%%%%%%%%%%%%%%%%%%%%%
\begin{figure}	\label{fig:mat_comp}
	\centering
	\setcounter{subfigure}{0}
	\begin{subfigure}{0.35\columnwidth}		\includegraphics[width=\linewidth, height = 0.7\linewidth]{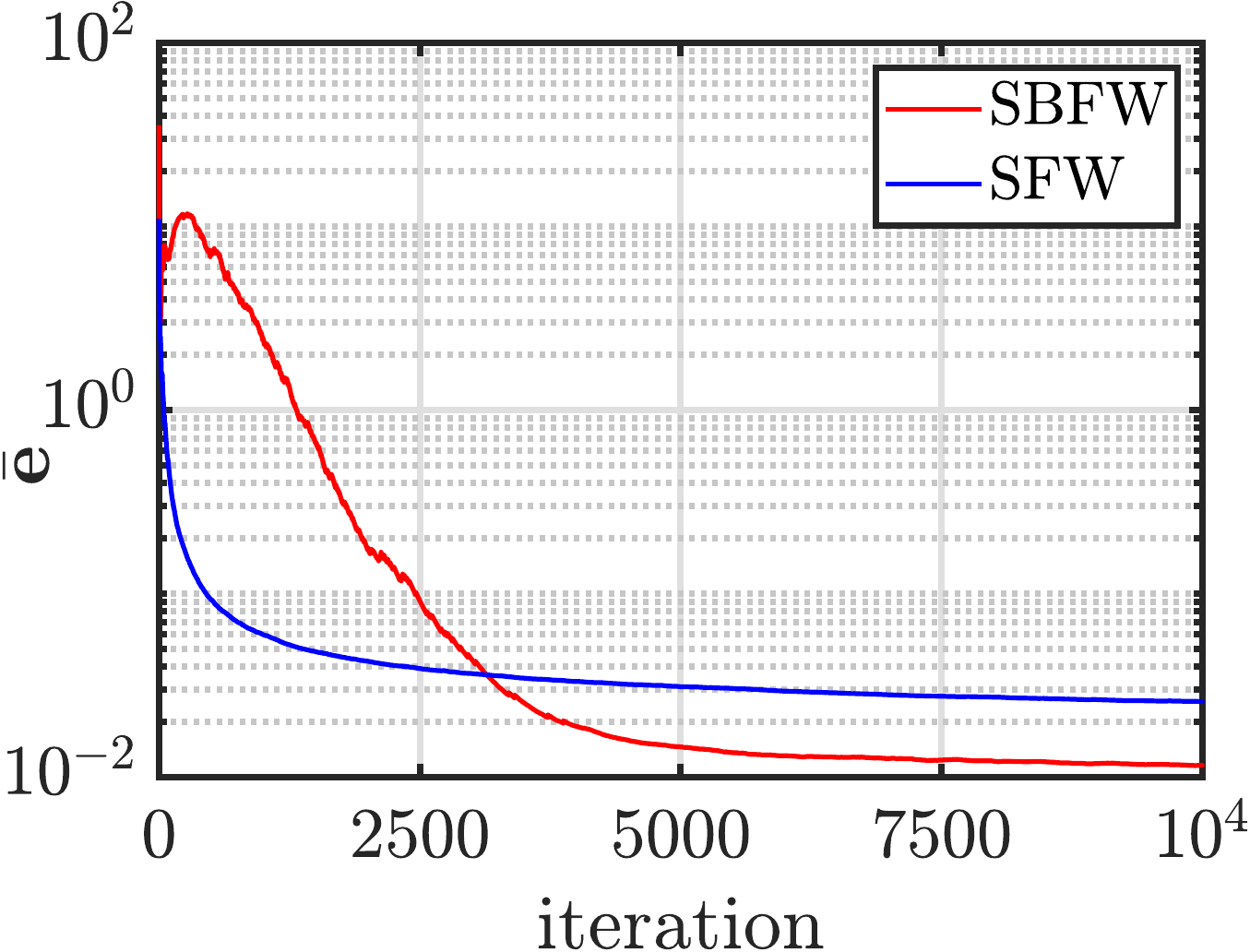}
		\caption{Normalized error.}
		\label{NE}
	\end{subfigure}
	\begin{subfigure}{0.35\columnwidth}		\includegraphics[width=\linewidth, height = 0.7\linewidth]{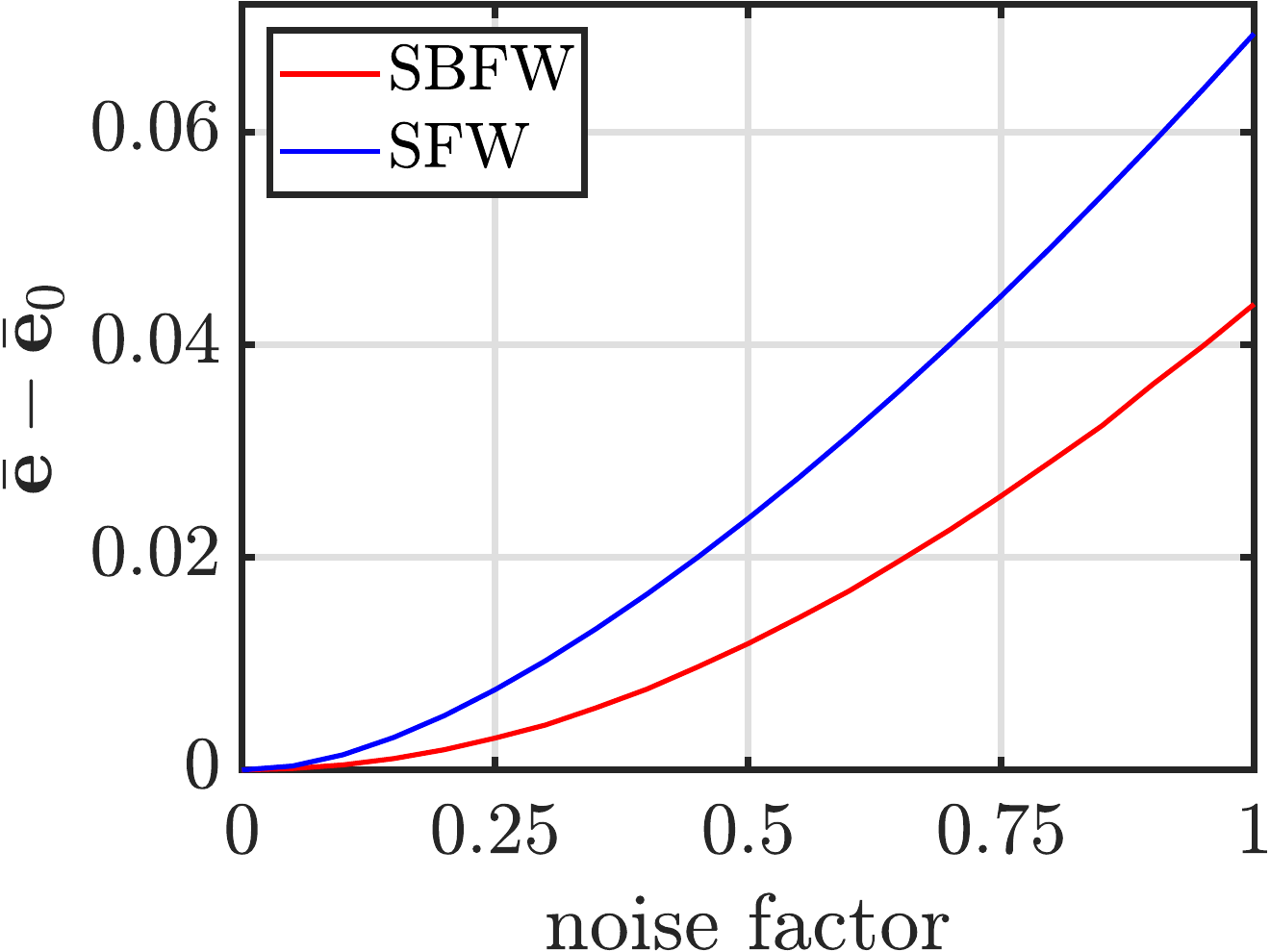}
		\caption{Error variation vs noise factor.}
		\label{noise}
	\end{subfigure}
	\caption{ This figure shows the benefit of bilevel formulation (solved via SBFW) for the matrix completion problem as compared to the standard single level problem (solved via SFW \cite{mokhtari2020stochastic}). We note that SBFW achieves a lower normalized error as compared to SWF (Fig. \ref{NE}). This advantage of SBFW is further confirmed in Fig. \ref{noise}, which plots the distance from the optimal value $\|\bar\e-\bar\e_0\|$ with respect to the noise factor in the matrix completion problem. The proposed algorithm is able to better handle the noise in the observed matrices by a significant margin. }
\end{figure}
\section{Numerical Experiments} In this section, first, we consider the  problem of low-rank matrix completion formulated in \eqref{mat_comp} to illustrate the performance of our proposed SBFW algorithm. Then we perform another experiment on sparse policy value evaluation to highlight the importance of the proposed SCFW algorithm.  
All the experiments are performed in MATLAB R2018a with Intel(R) Core(TM) i7-8550U CPU $@$ 1.80GHz.
\subsection{Importance of SBFW (cf. Algorithm \ref{alg:SBFW})}\label{sec:application}
%\subsection{Matrix Completion with Denoising}
{We start with emphasizing that there are two  challenges (C1 and C2) involved in solving the problem in \eqref{mat_comp};  C1: the problem is stochastic bi-level in nature, and C2: the nuclear norm constraints are expensive to project onto. We address both the challenges by solving  the problem in \eqref{mat_comp} via the proposed SBFW algorithm.  Note that \eqref{mat_comp} is of the form \eqref{first-level-sto} with $f(\cdot)=g(\cdot)=\|\cdot\|_F^2$, and $(\theta_t,\xi_t)$ being the  independent  random
	subset of entries revealed at every iteration $t$. }

\textbf{Synthetic Dataset Results:} We start with performing  experiments to  highlight the importance of denoising in matrix completion which lead to the bi-level formulation of the problem in \eqref{mat_comp}.  We follow the experimental setting of \cite{mokhtari2020stochastic}  and start with forming an observation matrix as $\M=\hat{\X}+\Eb$. Here, $\hat{\X}=\W\W^T$ with $\W\in \mathbb{R}^{n \times r}$ containing normally distributed independent entries, and the noise matrix $\Eb=\hat{n}(\L+\L^T)$ where $\L\in \mathbb{R}^{n \times n}$ contains normally distributed independent entries and $\hat{n}\in (0,1)$, is the noise factor. For the simulations, we set $n=250$, $r=10$, and $\alpha=\|\hat{\X}\|_*$. Further, we define the set of observed entries $\Omega$  by sampling $\M$ uniformly at random with probability $0.8$.

 We start with setting $\hat{n}=0.5$ and solve the problem in \eqref{mat_comp} using the proposed SBFW algorithm  and compare it with the state-of-the-art conditional gradient method called SFW \cite{mokhtari2020stochastic}. For SBFW, we set the step sizes as dictated  in  theory and set $\lambda_1=\lambda_2=0.05$, while for SFW, the step sizes are set as defined in \cite{mokhtari2020stochastic}. We use a batch size of $b=250$ for both the algorithms and run them for $10^4$ iterations. The performance is analyzed in terms of normalized error, defined as
\begin{align}
\bar{\textbf{e}}=\frac{\sum_{(i,j)\in \Omega}(\X_{i,j}-\hat{\X}_{i,j})^2}{\sum_{(i,j)\in \Omega}(\hat{\X}_{i,j})^2}
\end{align}
where $\X$ is the output generated by the algorithm. 

The evolution of the normalized error is shown in Fig.~\ref{NE}.
 Since SFW is a projection-free algorithm, it addresses the challenge C2 and converges fast but to a suboptimal point, while SBFW converges to better accuracy but slowly.  Hence, the proposed algorithm solves the bi-level matrix completion problem efficiently (addressed challenge C1). This also justifies the claim that additional denoising step improves the quality of matrix completion. Then, to further investigate the effectiveness of the proposed algorithm in addressing C1, we solve the problem for different noise factor $\hat{n}\in (0,1)$ and compare the normalized error in the solution at each $\hat{n}$ with the normalized error obtained for zero noise case (i.e. $\hat{n}=0$) denoted as $\bar{\textbf{e}}_0$ (shown in Fig.~\ref{noise}). We note that the growth in the error difference is much slower for SBFW than SFW, which shows the effectiveness of the proposed algorithm in dealing with noise.   
% \usepackage{multirow}

\begin{figure}[t]
	\centering
	\includegraphics[scale=0.6]{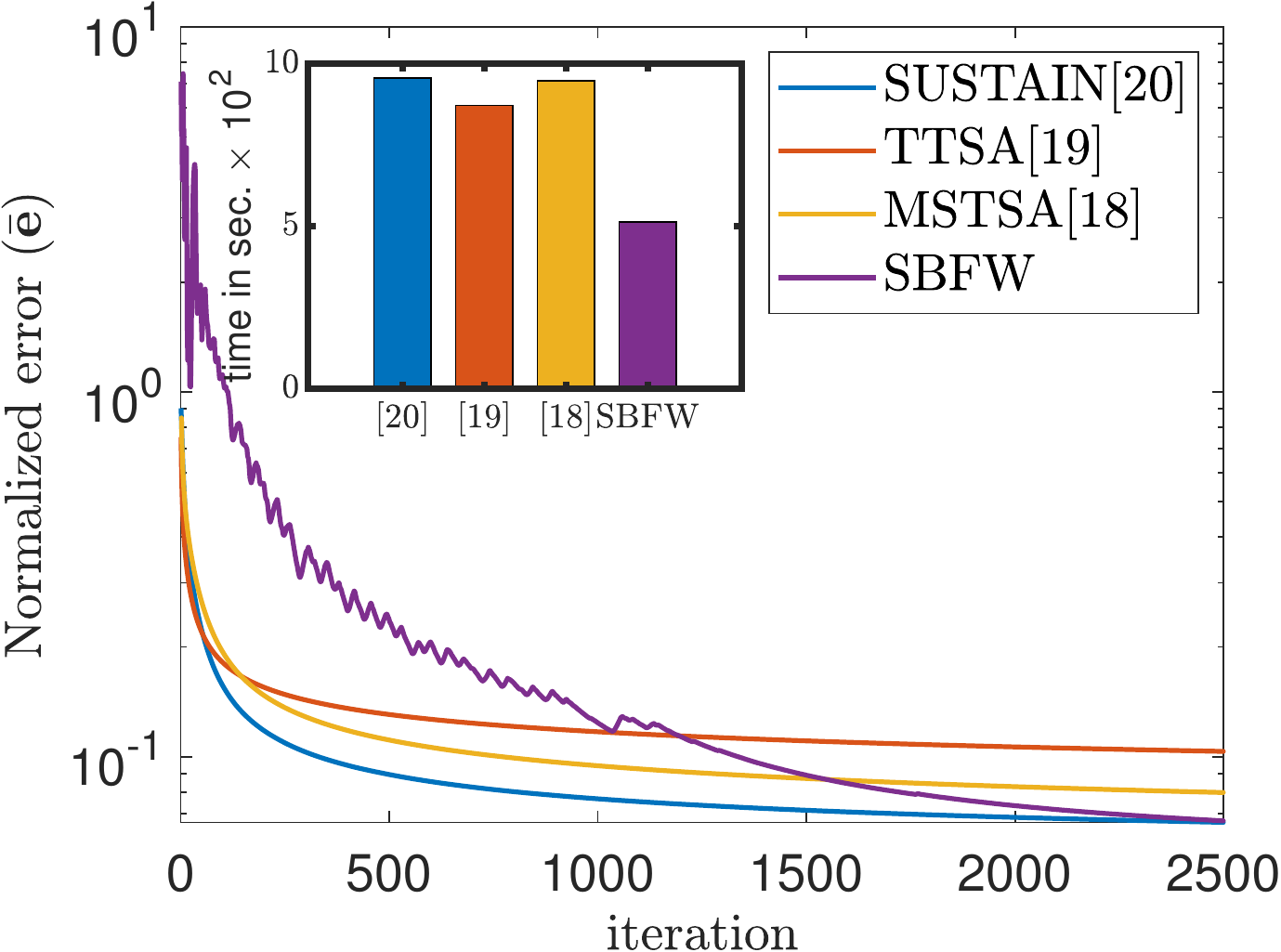}
	\caption{This figure compares the performance of the proposed SBFW algorithm  for matrix completion problem on MovieLens $100$k dataset with the other state-of-the-art algorithms such as TTSA \cite{hong2020two}, MSTSA \cite{khanduri2021momentum}, and SUSTAIN \cite{khanduri2021near}.We note that SBFW takes the least amount of time (shown in the inner bar plot) to attain same level of accuracy.}
	\label{fig:objvscv}
\end{figure}

\begin{figure}[t!]
	\centering
		\includegraphics[scale=0.6]{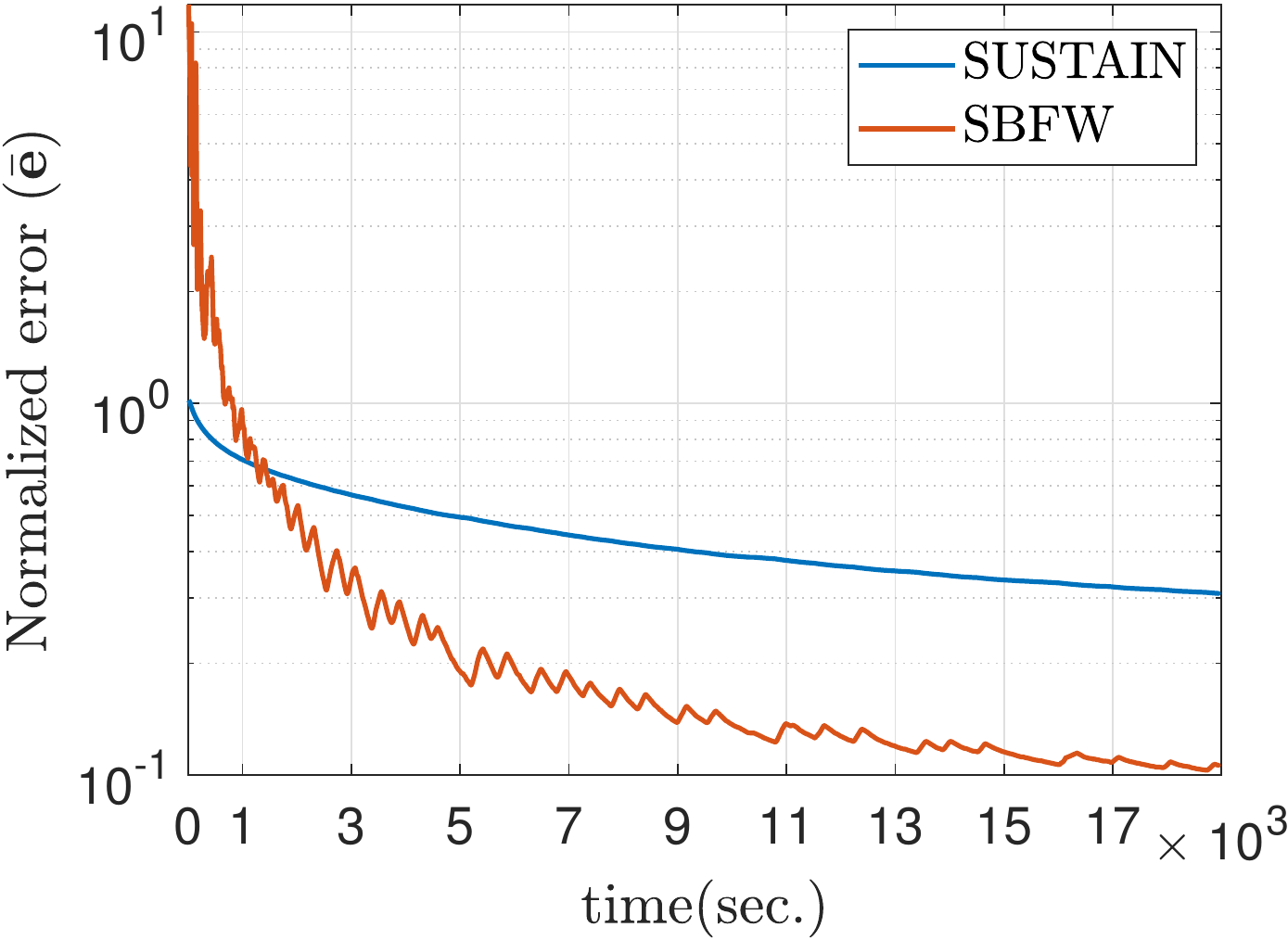}
	\caption{This figure compares the normalized error with respect to computation time required for  SBFW and SUSTAIN\cite{khanduri2021momentum} on MovieLens $1$M dataset. We note that eve thought SUSTAIN has a better theoretical convergence guarantees, it performs poorly in terms of computation time required to achieve a particular level of normalized error.  For instance, to achieve normalized error of $0.3$, SUSTAIN requires $31$ minutes while SBFGW requires only $5$ minutes. So there is a percentage improvement of $15\%$ by the proposed algorithm.  }
	\label{big_time}
\end{figure}

\textbf{Real Dataset Results:} To test the scalability of our proposed projection-free algorithm, we run experiment over large size matrices of  MovieLens\footnote{\href{https://grouplens.org/datasets/movielens/}{https://grouplens.org/datasets/movielens/}} datasets, which contains user ratings of movies ranging from $0$ to $5$. Since, SBFW is  a single loop algorithm, we compare its performance of SBFW with other state-of-the-art single loop projection-based bilevel algorithms such as SUSTAIN \cite{khanduri2021near}, TTSA \cite{hong2020two}, and MSTSA \cite{khanduri2021momentum}.  We start with Movielens $100$k dataset of $10^5$ ratings from $1000$ users for $1700$ movies. This dataset is denoted by observation matrix $\M$ of size $1000\times1700$. For the simulations, we define the set of observed entries $\Omega$  by sampling the matrix uniformly at random from $\M$ with a batch size of $b=5000$. Fig.~\ref{fig:objvscv} plots the evolution of normalized error for $2500$ iterations for all the algorithms. We note that the proposed algorithm is not the best in terms of the convergence rate when compared to projection-based schemes, which is expected from the slower theoretical convergence rates. However, when compared in terms of amount of clock time required to achieve the same level of normalized error (Fig.~\ref{fig:objvscv}),  he proposed scheme outperforms the other state-of-the-art methods. This gain is coming from the projection-free nature of the proposed algorithm, and we no longer required to perform a complicated projection at each iteration. 
To further highlight the importance of projection-free bilevel algorithm in practice, we perform additional experiments on larger dataset (of MovieLens $1$M) which contains $1$ million ratings from $6000$ users and for $4000$ movies. We plot the evolution of normalized error with time in Fig.~\ref{big_time}, where we only compare SBFW against SUSTAIN which is the state-of-the art projection based bilevel algorithm. It is interesting to note that even though SUSTAIN has a better theoretical convergence rate, it shows inferior performance in actual computation time (due to the projection operation) compared to SBFW as evident from Fig.~\ref{big_time}.  In Table \ref{time_table},  we provide  computation time comparisons (to complete $10^3$ iteration) of both the algorithms (under same settings) over different real datasets. Observe that for large data set, SBFW is approximately $10\times$ faster than the SUSTAIN and exhibit an improvement upto $82\%$ in the computation time. This performance gain in terms of computation time comes from the fact that other methods require to perform projections over nuclear norm at each iteration which is computationally expensive due to the computation of full singular value decomposition. In contrast, SBFW solves only a single linear program at each iteration, which only requires the computation of  singular vectors corresponding to the highest singular value. 

\begin{table}[h]
	\centering\scalebox{1}{
		\begin{tabular}{|c|c|c|c|cc|c|}
			\hline
			\multirow{2}{*}{Dataset} & \multirow{2}{*}{\#users} & \multirow{2}{*}{\#movies} & \multirow{2}{*}{\#ratings} & \multicolumn{3}{c|}{Time}           \\ \cline{5-7} 
			&                          &                           &                            & \multicolumn{1}{c|}{SUSTAIN} & SBFW & $\%$imp. \\ \hline
			Movielens 100k           & 1000                     & 1700                      & $10^5$      & \multicolumn{1}{c|}{554 sec.}     & 236 sec. & 57\%   \\ \hline
			Movielens latest         & 600                      & 9000                  & $10^5$      & \multicolumn{1}{c|}{66.6 mins.}     & 12.9 mins. & 81\%   \\ \hline
			Movielens 1M             & 6000                     & 4000                      & $10^6$      & \multicolumn{1}{c|}{10.16 hrs.}     & 1.82 hrs.  & 82\% \\ \hline
	\end{tabular}}
	\caption{Comparison of computation time of the proposed algorithm SBFW and the state-of-the-art projection based algorithm SUSTAIN over different size of real data sets.}
	\label{time_table}
\end{table}
We have analyzed the performance of SBFW on different datasets but a remark is due on the importance of proposed SCFW and studying its convergence rate results separately in this work. This is detailed next.
\subsection{Importance of SCFW (cf. Algorithm \ref{alg:SCFW})}\label{sparse_policyEXP}
In this subsection, we consider the problem of policy evaluation which lies at the heart of reinforcement learning as discussed in \cite{wang2017accelerating}. We perform this experiments to justify the claim that for compositional stochastic problems, the proposed SCFW (cf. Algorithm \ref{alg:SCFW}) exhibits a faster convergence rate as compared to SBFW (cf. Algorithm \ref{alg:SBFW}). 
%
%In this section, the proposed Algorithm \ref{alg:SCFW} (SCFW) is applied to the reinforcement learning application. We consider the problem of estimating the value function at a particular state often termed as policy learning. It is an important in-between step to find the optimal policy of any reinforcement learning problem. For instance, in Q-learning approach, the agent initially explores many possibilities to find the best action at a particular state such that the value function is maximized. However, evaluating value function is a dynamic decision problem solving which is not practical when the number of states is large. When the complete knowledge of transition probabilities is given, off-policy learning has been considered in \cite{white2016investigating} while the online version of policy learning has been discussed in \cite{wang2017accelerating}. We follow the similar approach as presented in \cite{wang2017accelerating}. 
%
%However, we focus on the specific experiment where the constraint make difficult apply projection or proximal steps.   
%
Consider a Markov Decision Process (MDP) with finite state space $\mathcal{S}$ and action space $\mathcal{A}$.  For a fixed policy $\pi$ which maps the current state $s$ to action $a\in\mathcal{A}$, the value function $V^{\pi}(s)$ for state $s$ is given by 
\begin{align}\label{valfun}
V^{\pi}(s) = \mbE\left\{r_{s,\hat{s}} + \gamma V^\pi(\hat{s})|s,\pi\right\} \forall s,\hat{s}\in\mathcal{S},
\end{align}
where $r_{s,\hat{s}}$ is the reward for transitioning from $s$ to $\hat{s}$, $\gamma\in (0,1)$ is the discount factor, and the expectation is taken over all the possible future states $\hat{s}$ conditioned on current state $s$ and policy $\pi$. Looking at the bellman equation in \eqref{valfun}, it is clear that evaluating the value function at all the states is impractical when $|\mathcal{S}|$ is moderately large. Hence, we consider a  linear function approximation for the value function such that $V^{\pi}(s) = \mathbf{\phi}_s^T\w^\star$, for some $\w^\star \in \mbR^m$. Here,  $\mathbf{\phi}_s \in \mbR^m$ denotes the $m$ dimensional state features for state $s$. The goal is to learn an optimal $\w^\star$ to obtain a suitable linear function approximation for the value function. The problem of finding $\w^\star$ can be formulated as 
\begin{align}\label{Prob}
\min_{\|\w\|_1 \leq \alpha} \sum_{s=1}^{S}\left(\mathbf{\phi}_s^T\w - q_{\pi,\hat{s}}(\w)\right)^2,
\end{align} 
where 
\begin{align*}
q_{\pi,\hat{s}}(\w) &= \mbE\left\{r_{s,\hat{s}} + \gamma \mathbf{\phi}_{\hat{s}}^T\w|s,\pi\right\} = \sum_{\hat{s}}P^\pi_{s\hat{s}}\left(r_{s,\hat{s}} + \gamma \mathbf{\phi}_{\hat{s}}^T\w\right).
\end{align*}
The constraint $\|\w\|_1 \leq \alpha$ in \eqref{Prob} for some sparsity tuning parameter $\alpha>0$ is useful in practice, where the sparsity needs to be ensured. For instance, when the number of states in $\mathcal{S}$ is large, the features of each state would become large, thus making the dimension of the optimization variable large. We propose to solve \eqref{Prob} in a  projection-free manner using the proposed SCFW algorithm. In the current context, the complete knowledge of  transition probabilities $\left\{P^\pi_{s\hat{s}}\right\}$ is not known but revealed sequentially.
%\begin{figure}[t!]
%	\centering
%	\setcounter{subfigure}{0}
%	\begin{subfigure}{0.49\columnwidth}		\includegraphics[width=\linewidth, height = 0.7\linewidth]{images/RL}
%		\label{resobj}
%	\end{subfigure}
%	\begin{subfigure}{0.49\columnwidth}		\includegraphics[width=\linewidth, height = 0.7\linewidth]{images/conVio}
%		\label{rescon}
%	\end{subfigure}
%		\caption{(a) Convergence results for ASCPG and SCFW algorithms, (b) Violation of constraint for different $\lambda$'s for ASCPG.}
%	\label{resrein}
%\end{figure}
\begin{figure}[t!]
	\centering
	\setcounter{subfigure}{0}
	\begin{subfigure}{0.5\columnwidth}		\includegraphics[width=\linewidth, height = 0.7\linewidth]{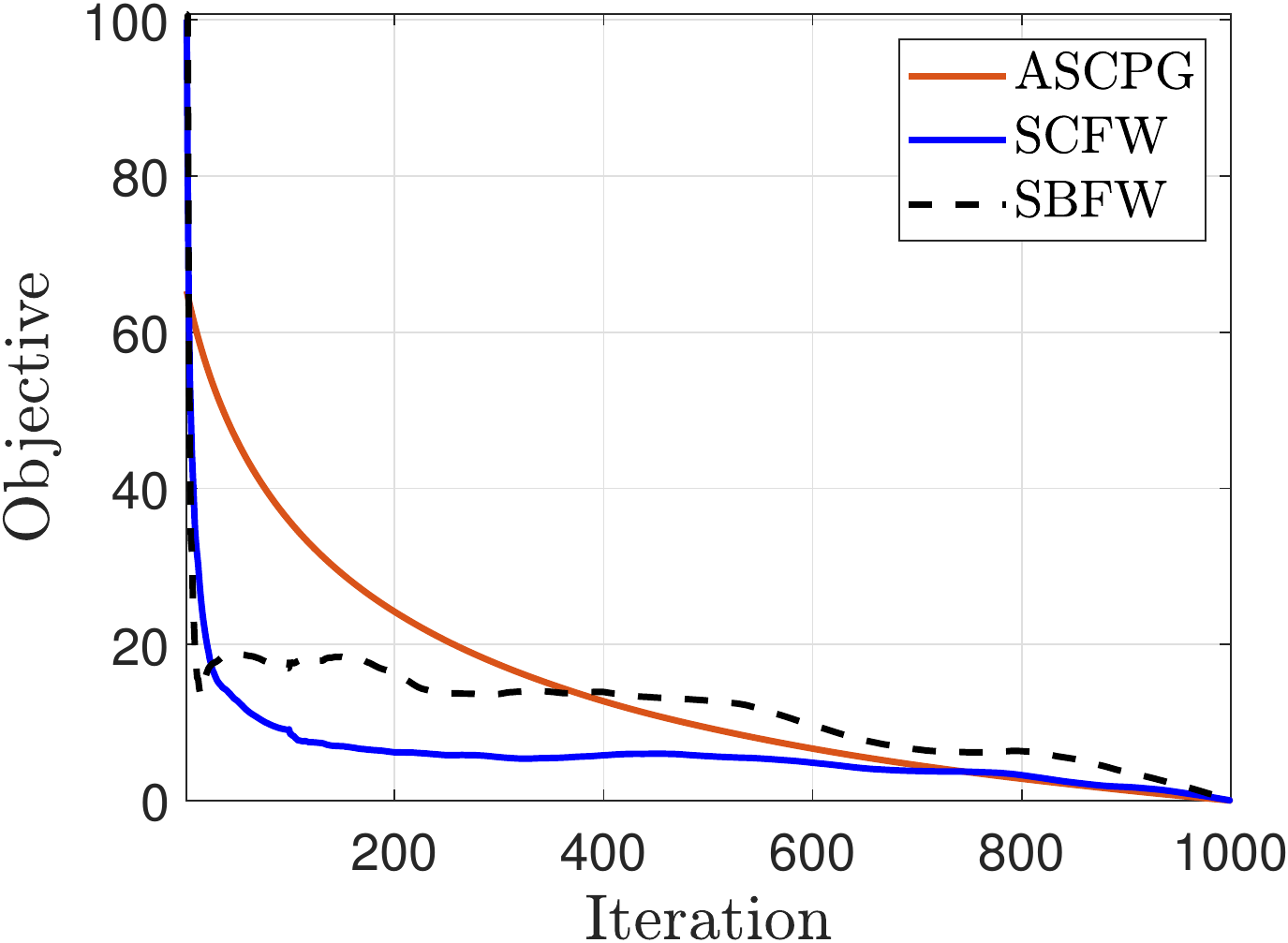}
		\label{resobj}
	\end{subfigure}
	\caption{Convergence results for ASCPG, SCFW and SBFW algorithms for sparse policy value evaluation (Sec.\ref{sparse_policyEXP}).}
	\label{resrein}
\end{figure}
For the experiments, we consider the number of states $|\mathcal{S}| = 100$ with $3$ actions available at every state. Given a pair of state and action, the agent can move any one of the next possible states. The transition probabilities and rewards for each transition are uniformly sampled in $[0, 1]$. For the purposeful behavior, we follow the same approach as in \cite{white2016investigating}, where the agent favors a single action at each state. Out of 3, one action is randomly selected and is assigned with the probability $0.9$, and the others are evenly assigned with probabilities. The feature vector of each state has dimension $m = 100$. However, the additional $\ell_1$ constraint ensures the sparsity in optimization variable $\w$.

 In our experiment, we set $\alpha=10^{-1}$ and run both the proposed SCFW and SBFW algorithm for 1000 iterations. The results of the experiment is reported  in terms of objective convergence $\|\w_t-\w^\star\|_2$ in Fig. \ref{resrein}. Here, $\w_t$ is the algorithmic solution at iteration $t$. The optimal solution $\w^\star$ is obtained by running the algorithm for $10^5$ iterations, given that  the complete information of transition probabilities is known. Observe that under the same setting, SCFW performs better as compared to SBFW. This performance gain is actually obtained due to the additional tracking of the inner function (see \eqref{inner_fucn_track}) in SCFW algorithm. We further provide comparison with the proximal ASC-PG algorithm \cite{wang2017accelerating}, where $\lambda\|\w\|_1$ is added as a  regularizer to the objective function. 
 The parameter $\lambda$ is tuned so that  $\|\w\|_1 \leq \alpha$ holds. As expected, all the algorithms have managed to reduce the error with iteration. However, note that ASC-PG is solving a relaxed version of the constraint problem by including the $\ell_1$ norm constraint as regularizer, that requires tuning of parameter $\lambda$. In fact, it might be possible to attain a better convergence plot for ASC-PG by tuning $\lambda$ further, but at an additional cost of constraint violation. In contrast,  the proposed algorithm is converging while strictly satisfying the constraint.
% To demonstrate this, we have run
%ASC-PG algorithm for different values of $\lambda$ (see Fig. 4b).
%Constraint violations are more when $\lambda$ is less. Hence, for
%the stricter constraint, more tuning is required. However, the
%proposed algorithm does not need that careful tuning and
%since it is not required to execute proximal step, the algorithm
%obtains the optimal value with lesser computational resources. 
Note that this experiment is just a representative example for SCFW and comparing with the gamut of RL algorithms is beyond the scope of the paper and hence we provide comparison only with the state-of-the-art proximal compositional algorithm \cite{wang2017accelerating}.

\section{Conclusion} \label{sec:conclusion}
This paper presents the first projection-free algorithm for stochastic bi-level optimization problems with a strongly convex inner objective function. We utilize the concept of momentum-based tracking to track the stochastic gradient estimate and establish the oracle complexities of the proposed SBFW algorithm for the convex and non-convex outer
objective functions. We also develop the first projection-free algorithm called SCFW for stochastic compositional problems.  We show that tracking both the inner function and the gradient of the objective function with momentum technique reduces  approximation noise which eventually helps obtain the optimal convergence rates. Numerical results show that the proposed projection-free variants have a significantly reduced wall-clock times as compared to their projection-based counterparts. 

\appendices
\section{Proof of Theorem \ref{TH:convex-bilevel}}\label{proof:theorem_convex_bilevel}
From the  initialization of variable $\x$, we have $\x_1\in \cX$. Also since we obtain $\s_t$ solving a linear minimization problem over the set $\cX$, we have $\s_t \in \cX$. Thus, $\x_{t+1}$ which is a convex combination of $\x_t$ and $\s_t$, i.e. $\x_{t+1}=(1-\eta_{t+1})\x_t+\eta_{t+1} \s_t$ will also lie in the set $\cX$. Hence $\x_{T+1}\in \cX$ and $\hat{\x}\in \cX$. Now, starting with  definition of $Q(\cdot)$, we have  $ Q(\x)= \mbE_{\theta}[f(\x,\y^\star(\x);\theta)]$. Also note that we have set $k=\frac{qL_g}{\mu_g}(\log(1+t))$, this ensures that the condition $\beta_t\leq \frac{C_{\x\y}C_{\y}}{\mu_g(t+1)^q}$ required in the analysis of Corollary \eqref{coroSBFW} is satisfied. Hence, we can use results from Corollary \eqref{coroSBFW} with $q=2/3$ for convex case and $q=1/2$ for non-convex case. 

\subsubsection{Proof of Statement (i) (Convex case)}
Using the smoothness assumption of $Q$ we can write
\begin{align}\label{smoothQ}
&Q(\x_{t+1})\!-\!Q(\x_t)\!\leq\! \ip{\n Q(\x_t),\x_{t+1}-\x_t}+\frac{L_Q}{2}\norm{\x_{t+1}\!-\!\x_t}^2\nonumber\\&=\eta_{t}\ip{\n Q(\x_t),\s_t-\x_t}+\frac{L_Q\eta_{t}^2}{2}\norm{\s_t-\x_t}^2,
\end{align}
where $L_Q=\frac{(L_{f_{\y}}+L)C_{\x\y}}{\mu_g}+L_{f_{\x}}+C_{\y}\left[\frac{L_{g_{\x\y}}C_{\y}}{\mu_g}+\frac{L_{g_{\y\y}}C_{\x\y}}{\mu_g^2}\right]$ (see Lemma \ref{Lemma:various_smoothness}). Here, in the last expression we have replace term $x_{t+1}-\x_t=\eta_{t}(\s_t-\x_t)$. Now adding and subtracting $\eta_{t}\ip{\d_t, \s_t-\x_t}$ in \eqref{smoothQ} we get
\begin{align}\label{smooth2Q}
Q(\x_{t+1})&\leq Q(\x_t)+\eta_{t}\ip{\n Q(\x_t)-\d_t,\s_t-\x_t}+\eta_{t}\ip{\d_t,\x^{\star}-\x_t}+\frac{L_Q\eta_{t}^2D^2}{2},
\end{align}
here in last the inequality is obtained using optimality of  $s_t$. Now introducing $\eta_{t}  \ip{\n Q(\x_t),\x^\star -\x_t}$ in RHS of \eqref{smooth2Q} and regrouping the terms we obtain
\begin{align}\label{smooth3Q}
Q(\x_{t+1})-\frac{L_Q\eta_{t}^2D^2}{2}&\leq Q(\x_t)+\eta_{t}\ip{\n Q(\x_t)-\d_t,\s_t-\x^\star}+\eta_{t}\ip{\n Q(\x_t),\x^\star-\x_t}\nonumber\\&\leq Q(\x_t)+\eta_{t} D\norm{\n Q(\x_t)-\d_t}+\eta_{t}\ip{\n Q(\x_t),\x^{\star}-\x_t}\nonumber\\&\leq Q(\x_t)+\eta_{t} D\norm{\n Q(\x_t)-\d_t}-\eta_{t}(Q(\x_t)-Q(\x^\star)),
\end{align}
here in the second inequality we use bound $\eta_{t}\ip{\n Q(\x_t)-\d_t,\s_t-\x^\star}\leq \eta_{t}\norm{\n Q(\x_t)-\d_t} \norm{\s_t-\x^\star}\leq \eta_{t} D\norm{\n Q(\x_t)-\d_t}$ and in last inequality we used the bound $\ip{\n Q(\x_t),\x^{\star}-\x_t}\leq Q(\x^\star)-Q(\x_t)$.
Subtracting $Q(\x^\star)$, taking expectation and using $\mbE \norm{X} \leq \sqrt{\mbE\norm{X}^2}$ we get
\begin{align}\label{Lemma for Q final eq2}
\mbE[Q(\x_{t+1})&-Q(\x^\star)]
\leq (1-\eta_{t})\mbE[Q(\x_t)-Q(\x^\star)]+ \eta_{t} D\sqrt{\mbE\norm{\n Q(\x_t)-\d_t}^2}+\frac{L_Q\eta_{t}^2D^2}{2}.
\end{align}
Further, setting $q=2/3$ hence, $\eta_t=\frac{2}{t+1}$ and using Corollary \ref{coroSBFW}, we can bound the second term of \eqref{Lemma for Q final eq2} $
\eta_tD\sqrt{\mbE\norm{\n Q(\x_t)-\d_t}^2} \leq \frac{2D\sqrt{C_1 }}{(t+1)^{4/3}}.
$
which gives
\begin{align}\label{Lemma for Q final eq3}
\mbE[Q(\x_{t+1})-Q(\x^\star)]&
\leq \left(1-\frac{2}{t+1}\right)\mbE[Q(\x_t)-Q(\x^\star)]+ \frac{2 D\sqrt{C_1}}{(t+1)^{\frac{4}{3}}}+\frac{2L_QD^2}{(t+1)^{2}}.
\end{align}
Multiplying both side by $t(t+1)$ we can write
\begin{align}\label{Qcc }
t(t+1)\mbE[Q(\x_{t+1})-Q(\x^\star)]&
\leq t(t-1)\mbE[Q(\x_t)-Q(\x^\star)]+ \frac{2t D\sqrt{C_1}}{(t+1)^{\frac{1}{3}}}+\frac{2tL_QD^2}{t+1}\nonumber\\
&\leq t(t-1)\mbE[Q(\x_t)-Q(\x^\star)]+ 2 D\sqrt{C_1}(t+1)^{\frac{2}{3}}+2L_QD^2,
\end{align}
Summing for $t=1,2,\cdots,T$ and rearranging we get
\begin{align}\label{Qcc1}
\mbE[Q(\x_{T+1})-Q(\x^\star)]
&\leq \frac{1}{T(T+1)}\left(\frac{6}{5} D\sqrt{C_1}(T+1)^{\frac{5}{3}}+2L_QD^2T\right)\nonumber\\
&\leq \frac{12 D\sqrt{C_1}}{5(T+1)^{\frac{1}{3}}}+\frac{2L_QD^2}{(T+1)},
\end{align}
here we use the fact that $\sum_{t=1}^{T}(t+1)^{2/3}\leq \frac{3}{5}(T+1)^{5/3}$.
\subsubsection{Proof of Statement (ii) (Non-convex case)}
Again starting with the smoothness assumption of $Q$ we can write
\begin{align}\label{smooth non convexQ}
Q(\x_{t+1})-Q(\x_t)&\leq \ip{\n Q(\x_t),\x_{t+1}\!-\!\x_t}\!+\!\frac{L_Q}{2}\norm{\x_{t+1}\!-\!\x_t}^2\nonumber\\&=\eta_{t}\ip{\n Q(\x_t),\s_t-\x_t}+\frac{L_Q\eta_{t}^2}{2}\norm{\s_t-\x_t}^2\\&\leq \eta_{t}\ip{\n Q(\x_t),\s_t-\x_t}+\frac{L_Q\eta_{t}^2D^2}{2}\nonumber\\&= \eta_{t}\ip{\n Q(\x_t)-\d_t,\s_t-\x_t}+\eta_{t}\ip{\d_t,\s_t-\x_t}+\frac{L_Q\eta_{t}^2D^2}{2},\nonumber
\end{align}
where in the second expression we have replace term $x_{t+1}-\x_t=\eta_{t}(s_t-\x_t)$ and in third inequality we used compactness assumption of set while in the last expression we introduced $\eta_{t}\ip{\d_t, s_t-\x_t}$. Next,  introducing the following quantity
$\hat{\upsilon}_t=\arg\max_{\upsilon  \in \mathcal{X}} \ip{\upsilon-\x_t,-\n Q(\x_t)},
$ and  using the optimality of $s_t$ that is $s_t=\arg \min_{\s\in \mathcal{X}}\ip{\d_t,\s}$, we write  \eqref{smooth non convexQ} as
\begin{align}\label{smooth non convex2Q}
Q(\x_{t+1})- Q(\x_t)&\leq\eta_{t}\ip{\n Q(\x_t)-\d_t,\s_t-\x_t}+\eta_{t}\ip{\d_t,\hat{\upsilon}_t-\x_t}+\frac{L_Q\eta_{t}^2D^2}{2}\nonumber\\&= -\eta_{t}\mathcal{G}(\x_t)+Q(\x_t)+\eta_{t}\ip{\n Q(\x_t)-\d_t,\s_t-\x_t}\nonumber\\&\quad+\eta_{t}\ip{\d_t-\n Q(\x_t),\hat{\upsilon}_t-\x_t}+\frac{L_Q\eta_{t}^2D^2}{2}
\nonumber\\&\leq  -\eta_{t}\mathcal{G}(\x_t)+Q(\x_t)+2\eta_{t}D\norm{\d_t-\n Q(\x_t)}+\frac{L_Q\eta_{t}^2D^2}{2},
\end{align}
where the second expression is obtained by adding and subtracting the term $\eta_{t}\ip{\hat{\upsilon}_t-\x_t,-\n Q(\x_t)}$ while the last inequality is obtained using compactness assumption.

Rearranging  \eqref{smooth non convex2Q}, summing for $t=1,2,\cdots,T$ and taking expectation, we get
\begin{align}\label{sumQ}
&\sum_{t=1}^{T}\eta_{t}\mbE[\mathcal{G}(\x_t)]-\frac{L_QD^2}{2}\sum_{t=1}^{T}\eta_{t}^2\\&\leq Q(\x_1)-\mbE[Q(\x_{T+1})]+2D\sum_{t=1}^{T}\eta_{t}\mbE\norm{\d_t-\n Q(\x_t)}\nonumber\\&
\leq Q(\x_1)-Q(\x^{\star})+2D\sum_{t=1}^{T}\eta_{t}\mbE\norm{\d_t-\n Q(\x_t)},
\end{align}
here the last inequality follows from optimality of $\x^\star$. Using bound from  Corollary \ref{coroSBFW} for $q=1/2$, and  Jensen's equality we can write
$
\sqrt{\mbE\norm{\d_t-\n Q(\x_t)}^2}\leq \frac{\sqrt{C_1}}{(t+1)^{\frac{1}{4}}}.
$
%Note that our selection of $\eta_t=\eta=\frac{2}{(T+2)^{2/3}}$ satisfies the condition $\eta_t\leq\frac{2}{(t+2)^{2/3}}$ required in Lemma \ref{Lemma: ySBQW} and Lemma \ref{Lemma: f-d}. 
Now, setting $\eta_t=\frac{2}{(T+1)^{3/4}}$ in \eqref{sumQ} and rearranging we can write
\begin{align}
\mbE[\mathcal{G}(\hat{\x})]&\leq \frac{1}{T}\sum_{t=1}^{T}\mbE[\mathcal{G}(\x_t)]
\nonumber\\&\leq\frac{Q(\x_1)-Q(\x^\star)}{2T(T+1)^{-3/4}}+\sum_{t=1}^{T}\frac{2D\sqrt{C_1}}{T(t+1)^{1/4}}+\sum_{t=1}^{T}\frac{L_QD^2}{T(T+1)^{3/4}}\nonumber\\
&\leq \frac{Q(\x_1)-Q(\x^\star)}{2T(T+1)^{-3/4}}+\frac{8 D\sqrt{C_1}}{3T(T+1)^{-3/4}}+\frac{L_QD^2}{(T+1)^{3/4}}\nonumber\\
&\leq \frac{Q(\x_1)-Q(\x^\star)}{(T+1)^{1/4}}+\frac{16 D\sqrt{C_1}}{3(T+1)^{1/4}}+\frac{L_QD^2}{(T+1)^{3/4}}\nonumber\\
&:=\mathcal{O}((T+1)^{-1/4}),
\end{align}
here in the second inequality we use the fact that $\sum_{t=1}^{T}(t+1)^{-1/4}\leq \frac{4}{3}(T+1)^{3/4}$. 
\subsection{Proof of Lemma \ref{Lemma:ySCFW} }\label{proof:Lemma:ySCFW}
Setting $\Psi(\cdot)=h(\cdot)$ in Lemma \eqref{Lemma:general_tracking}, we can write
%We start with introducing $(1-\delta_t)h(x_{t-1})$ in the update equation \eqref{inner_fucn_track} and then take the norm square, which gives
%\begin{align}
%\norm{\y_{t+1}-h(\x_t)}^2=\|(1-\delta_t)(\y_t-h(\x_{t-1}))-(1-\delta_t)(\hht(\x_{t-1})-h(\x_{t-1}))+\hht(\x_t)-h(\x_t)\|^2.
%\end{align}
%Now, expanding the square, taking conditional expectation $\mbE_t=\mbE[(\cdot)|\mathcal{F}_t]$ and using the fact that $\mbE_t[h(\x_{t-1})-\hht(\x_{t-1})]=0$ and $\mbE_t[h(\x_{t})-\hht(\x_{t})]=0$ we obtain
\begin{align}\label{half y}
\mbE_t[\norm{\y_{t}-h(\x_{t})}^2]&\leq (1-\delta_t)^2\norm{(\y_{t-1}-h(\x_{t-1}))}^2+2(1-\delta_t)^2\mbE_t[\norm{\hht(\x_t)-\hht(\x_{t-1})}^2]\nonumber\\&+2\delta_t^2\mbE_t[\norm{\hht(\x_t)-h(\x_t)}^2].
\end{align} 
Now, taking total expectation of \eqref{half y} we get
\begin{align}\label{temp}
&\mbE\norm{\y_{t}-h(\x_t)}^2-(1-\delta_t)^2\mbE\norm{(\y_{t-1}-h(\x_{t-1}))}^2\nonumber\\&\leq 2\delta_t^2\sigma_h^2+2(1-\delta_t)^2\mbE\norm{\hht(\x_t)-\hht(\x_{t-1})}^2\nonumber\\&\leq 2\delta_t^2\sigma_h^2+2(1-\delta_t)^2M_h\mbE\norm{\x_t-\x_{t-1}}^2\nonumber\\&= 2\delta_t^2\sigma_h^2+2(1-\delta_t)^2M_h\eta_{t-1}^2\mbE\norm{\s_{t-1}-\x_{t-1}}^2
\nonumber\\&\leq  2\delta_t^2\sigma_h^2+2(1-\delta_t)^2\eta_{t-1}^2M_hD^2,
\end{align} 
here in the first inequality we have used Assumption \ref{ass:SCFW}(\ref{ass:inner_var}). In the second inequality we have used Assumption \ref{ass:SCFW}(\ref{ass:bounded_moment}) while the next expression follows from the update step $\x_{t+1}=(1-\eta_{t})\x_t+\eta_{t} \s_t$. The last inequality obtained using compactness assumption of set $\cX$.

\bibliographystyle{IEEEtran}
\bibliography{ref}

%
%                            

%%                        

%%%%%%%%%%%%%%%%%%%%%%%%%%%%
%%%%%%%%%%%%%%%%%%%%%%%%%%%%%%%%
%Supplementary
%%%%%%%%%%%%%%%%%%%%
%%%%%%%%%%%%%%%%%%%%

\newpage
\section*{Supplementary Material}

\section{Preliminaries}\label{key}

%\subsection{Notations} 
%
%First, we defined the compact notations we utilize in the convergence proofs. We denote column vectors with lowercase boldface $\x$, its transpose as $\x^\top$, and its Euclidean norm by $\norm{\x}$.  We use $\mbE_t:=\mbE[\cdot|\mathcal{F}_t]$ to denote
%the conditional expectation  with respect to given sigma field $\mathcal{F}_t$
%which contains all algorithm history (randomness) till step $t$.   
%% We use $\n$ operator without any subscript to represent gradient with respect to $\x$.
%%  
%

\subsection{Existing Results}
We start the discussion by mentioning some of the existing results in Lemma \ref{Lemma:various_smoothness} and \ref{Lemma:borw_result2} which are useful for the analysis in this paper.
\begin{lemma}\label{Lemma:various_smoothness}
	[\cite{ghadimi2018approximation}, Lemma 2.2]
	Under Assumption\ref{ass:SBFW}, following statements hold.
	\begin{enumerate}[(a)]
		\item For any $\x\in \cX$ and $\y \in \mathbb{R}^n$, 	\begin{align}\label{various_smoothness}
		\norm{\n S(\x,\y)-\n Q(\x)}\leq L\norm{\y^\star(\x)-\y},
		\end{align}
		where $L:=L_{f_{\x}}+\frac{L_{f_{\y}}C_{\x\y}}{\mu_g}+C_{\y}\left[\frac{L_{g_{\x\y}}}{\mu_g}+\frac{L_{g_{\y\y}}C_{\x\y}}{\mu_g^2}\right]$ and all the constants are as defined in  Assumption\ref{ass:SBFW}.
		
		\item The inner optimal solution $\y^\star(\x)$ is $\frac{C_{\x\y}}{\mu_g}$-Lipschitz continuous in $\x$, which implies that for any $\x_1,\x_2\in \cX$, it holds that	$\norm{\y^\star(\x_1)-\y^\star(\x_2)}\leq \frac{C_{\x\y}}{\mu_g}\norm{\x_1-\x_2}$.
		
		\item The gradient of outer objective $\n Q$ is $L_Q$-Lipschitz continuous in $\x$, which implies that for any $\x_1,\x_2\in \cX$, it holds that $\norm{Q(\x_1)-Q(\x_2)}\leq L_Q\norm{\x_2-\x_1}$ where $L_Q:=\frac{(L_{f_y}+L)C_{xy}}{\mu_g}+L_{f_x}+C_y\left[\frac{L_{g_{xy}}C_y}{\mu_g}+\frac{L_{g_{yy}}C_xy}{\mu_g^2}\right]$.
		
	\end{enumerate}
	
\end{lemma}

\begin{lemma}\label{Lemma:borw_result2}
	[\cite{khanduri2021near}, Lemma 4.1]
	Suppose Assumption\ref{ass:SBFW} holds, and the gradient estimate $h(\x,\y;\theta,\xi)$ is constructed with $k$ number of samples using \eqref{biased_estimate},  then\\
\noindent	(a) for any $\x\in \cX$ and $\y_1,\y_2 \in \mathbb{R}^n$, we have	
		\begin{align}
		\mbE_t\|h(\x,\y_{1};\!\theta_t,\!\xi_t)\!-\!h(\x,\y_{2};\!\theta_t,\!\xi_t)\|^2\!\leq\! L_k\mbE_t\|\y_{1}\!-\!\y_{2}\|^2
		\end{align}
		
		\noindent (b) for any  $\y \in \mathbb{R}^n$  and $\x_1,\x_2\in \cX$, we have	
		\begin{align}
		\mbE_t\|h(\x_{1},\y;\!\theta_t,\!\xi_t)\!-\!h(\x_2,\y;\!\theta_t,\!\xi_t)\|^2\!\leq L_k\mbE_t\|\x_{1}\!-\!\x_{2}\|^2
		\end{align}		
	where 
	$$L_k=2L^2_{f_{\x}}+\tfrac{6k[(L_g-\mu_g)^2(C^2_{g_{\x\y}}L^2_{f_{\y}}+C^2_{f_{\y}}L^2_{g_{\x\y}})+k^2C^2_{g_{\x\y}}C^2_{f_{\y}}L^2_{g_{\y\y}}]}{\mu_g(2L_g-\mu_g)}.$$
\end{lemma}
%The first part of Lemma \ref{Lemma:borw_result2} shows that $\mbE[h(\x_t,\y_t;\theta_t,\xi_t)]=\n S(\x_t,\y_t)+B_t$, with
%\begin{align}\label{bt}
%\norm{B_t}\leq \beta_t=\mathcal{O}\left(\left(1-\frac{\mu_g}{L_g}\right)^k\right),
%\end{align}
%while part (b) and (c) establishes the Lipschitzness of stochastic gradient estimate.
\subsection{General Inequalities}
Before proceeding towards the main analysis of this work, we first present and establish a general mathematical inequality in the form of Lemma \ref{Lemma:gen_seq} which will be useful to the analysis in this work. Further, we will present a general upper bound in Lemma \ref{Lemma:general_tracking} on the expected estimation error when the momentum-based method is employed to track the function or gradient.

\begin{lemma}\label{Lemma:gen_seq}
	Let $\psi_t$ be a sequence of real numbers which satisfy 
	\begin{align}\label{psi_claim}
	\psi_{t+1}=\left(1-\frac{c_1}{(t+t_0)^{r_1}}\right)\psi_{t}+\frac{c_2}{(t+t_0)^{r_2}}
	\end{align}
	for some $r_1\in (0,1]$ such that $r_1\leq r_2\leq 2r_1$, $c_1>1$, and $c_2\geq 0$. Then, $\psi_{t+1}$ would converge to zero at the following rate
	\begin{align}\label{gen_seq_rate}
	\psi_{t+1}\leq \frac{c}{(t+t_0+1)^{r_2-r_1}},
	\end{align}
	where $c=\max\{\psi_1 (t_0+1)^{r_2-r_1},\frac{c_2}{c_1-1}\}$.
\end{lemma}
\noindent\textbf{Proof:}
We prove Lemma \ref{Lemma:gen_seq} by induction. The base step of induction holds (for $t=0$) from the definition of $c$ which implies that $c\geq \psi_1 (t_0+1)^{r_2-r_1}$. Next, we assume that \eqref{gen_seq_rate} holds for $t=k$, which means
\begin{align}\label{gen_seq_rate_k}
\psi_{k+1}\leq \frac{c}{(k+t_0+1)^{r_2-r_1}}.
\end{align}
Now it remains to show that \eqref{gen_seq_rate} also holds for $t=k+1$. To proceed, we set $t=k+1$ in \eqref{psi_claim} to obtain
\begin{align}\label{psi_claim_k}
\psi_{k+2}=\left(1-\tfrac{c_1}{(k+t_0+1)^{r_1}}\right)\psi_{k+1}+\tfrac{c_2}{(k+t_0+1)^{r_2}}.
\end{align}
From the definition of $c$, it holds that $c_2\leq c(c_1-1)$ and we utilize the upper bound on $\psi_{k+1}$ from \eqref{gen_seq_rate_k} into \eqref{psi_claim_k} to obtain
\begin{align}\label{psi_claim_k2}
\psi_{k+2}\leq \left(1-\tfrac{c_1}{(k+t_0+1)^{r_1}}\right)\tfrac{c}{(k+t_0+1)^{r_2-r_1}}+\tfrac{c(c_1-1)}{(k+t_0+1)^{r_2}}.
\end{align}
Simplify and rearrange the terms in \eqref{psi_claim_k2} to write 
\begin{align}
\psi_{k+2}&\leq c\left(\frac{1}{(k+t_0+1)^{r_2-r_1}}-\frac{1}{(k+t_0+1)^{r_2}}\right).
\end{align}
Next, note that for general $p,\ q$, it holds that $\frac{1}{x^{p-q}}-\frac{1}{x^p}\leq \frac{1}{(x+1)^{p-q}}$, which allows us to write 
\begin{align}
\psi_{k+2}\leq \frac{c}{(k+t_0+2)^{r_2-r_1}}.
\end{align}
Thus, it holds that the statement in \eqref{psi_claim} holds for all $t\geq 0$.
%Next, we establish a special case of above result in Lemma \ref{Lemma_simplelemma}. \colr{In what sense is this special?} 
%%
%\begin{lemma}\label{Lemma_simplelemma}
% For some constant $A>0$ and $a\in[1,2]$, let $\{\psi_k\}_{k=1}^t$ be a sequence of numbers satisfying the following recursion
% %
%	\begin{align}
%	\psi_k &\leq \left(1 - \frac{1}{k}\right)^2\psi_{k-1} + \frac{A}{k^a} \label{simplerec1}.
%	\end{align}
%%
% Then it holds that $\psi_k \leq 4A/(k+1)^{a-1}$ for all $t \geq 1$. 
%\end{lemma}
%\begin{IEEEproof}
%Stat by	multiplying both sides of \eqref{simplerec1} by $k^2$ to obtain 
%	\begin{align}
%		k^2 \psi_k &\leq (k-1)^2\psi_{k-1} + Ak^{2-a} ,
%	\end{align}
%which holds for  $k \geq 1$. Take telescopic sum over $\tau = 1, \ldots, k$, we obtain
%	\begin{align}
%		\psi_k &\leq \frac{A}{k^2}\sum_{\tau=1}^k \tau^{2-a} 
%\leq A\frac{(k+1)^{3-a}}{(3-a)k^2}		\nonumber
%\\
%&\leq A\frac{(k+1)^2(k+1)^{3-a}}{(k+1)^2k^2}  
%		 \leq 4A\frac{(k+1)^{3-a}}{(k+1)^2} 		\nonumber
%		 \\
%		 &\leq \frac{4A}{(k+1)^{a-1}},
%	\end{align}
%	where we have used the inequality $\sum_{\tau=1}^{k} \tau^{n} \leq \frac{(k+1)^{n+1}-1}{(n+1)}$ for $n \geq 0$.
%	\end{IEEEproof}
\begin{lemma}\label{Lemma:general_tracking}
	Let us estimate function $\Psi(\bbx)=\mbE_{\xi}[\Psi(\bbx,\xi)]$ by $\y_t$ using step size $\delta_t$ as follows
	\begin{align}\label{gen_track}
	\y_{t}&=(1-\delta_t)(\y_{t-1}- \Psi(\x_{t-1},\xi_t))+ \Psi(\x_t,\xi_t).
	\end{align}
	Then the expected tracking error $\mbE_t[\norm{\y_{t}-\Psi(\x_t)}^2]$ satisfies 
	\begin{align}\label{gen0}
	\mbE_t[\norm{\y_{t}-\Psi(\x_t)}^2]&\leq(1-\delta_t)^2\norm{(\y_{t-1}-\Psi(\x_{t-1}))}^2\nonumber
	\\
	&\quad+2(1-\delta_t)^2\mbE_t[\norm{\Psi(\x_{t},\xi_t)-\Psi(\x_{t-1},\xi_t)}^2]+2\delta_t^2\mbE_t[\norm{\Psi(\x_{t},\xi_t)-\Psi(\x_t)}^2].
	\end{align} 
\end{lemma}
\textbf{Proof:} Consider the update equation in  \eqref{gen_track}, add/subtract the term $(1-\delta_t)\Psi(\x_{t-1})$ in the right hand side of \eqref{gen_track} to obtain
\begin{align}\label{gen_track2}
\y_{t}&=(1-\delta_t)(\y_{t-1}- \Psi(\x_{t-1},\xi))+ \Psi(\x_t,\xi_t)+(1-\delta_t)\Psi(\x_{t-1})-(1-\delta_t)\Psi(\x_{t-1}).
\end{align}
%s
Subtract $\Psi(\x_t)$ from both sides in \eqref{gen_track2} and take norm square:
\begin{align}
&\norm{\y_{t}-\Psi(\x_t)}^2=\|(1-\delta_t)(\y_{t-1}-\Psi(\x_{t-1}))\\&\quad-(1-\delta_t)(\Psi(\x_{t-1},\xi_t)-\Psi(\x_{t-1}))+\Psi(\x_{t},\xi_t)-\Psi(\x_t)\|^2.\nonumber
\end{align}
Now, expand the square and calculate conditional expectation $\mbE_t=\mbE[(\cdot)|\mathcal{F}_t]$ to obtain
\begin{align}
\mbE_t[\norm{\y_{t}-\Psi(\x_t)}^2]
&=(1-\delta_t)^2\norm{(\y_{t-1}-\Psi(\x_{t-1}))}^2
\\
&-2\langle(1-\delta_t)(\y_{t-1}-\Psi(\x_{t-1})),(1-\delta_t)(\mbE_t[\Psi(\x_{t-1})-\Psi(\x_{t-1},\xi_t)])+\mbE_t[\Psi(\x_{t})-\Psi(\x_{t},\xi_t)]\rangle
\nonumber
\\
&\!+\!\mbE_t\|(1\!-\!\delta_t)(\Psi(\x_{t-1},\xi_t)\!-\!\Psi(\x_{t-1}))\!+\!\Psi(\x_t)\!-\!\Psi(\x_{t},\xi_t)\|^2.\nonumber
\end{align}
Note that $\mbE_t[\Psi(\x_{t-1})-\Psi(\x_{t-1},\xi_t)]=0$ and $\mbE_t[\Psi(\x_{t})-\Psi(\x_{t},\xi_t)]=0$, which implies that
\begin{align}\label{gen1}
\mbE_t[\norm{\y_{t}-\Psi(\x_t)}^2]&=(1-\delta_t)^2\norm{(\y_{t-1}-\Psi(\x_{t-1}))}^2
+\!\mbE_t\norm{(1\!-\!\delta_t)(\Psi(\x_{t-1},\!\xi_t)\!-\!\Psi(\x_{t-1}))\!+\!\Psi(\x_t)\!-\!\Psi(\x_{t},\xi_t)}^2\nonumber\\&\leq(1-\delta_t)^2\norm{(\y_{t-1}-\Psi(\x_{t-1}))}^2+ 2(1-\delta_t)^2\mbE_t[\norm{\Psi(\x_{t},\xi_t)-\Psi(\x_{t-1},\xi_t)}^2]\nonumber\\&\quad+2\delta_t^2\mbE_t[\norm{\Psi(\x_{t},\xi_t)-\Psi(\x_t)}^2]
\end{align} 
where the last inequality holds due to the fact that $\mbE\norm{X-\mbE[X]+Y}^2\leq 2\mbE\norm{X}^2+2\mbE\norm{Y}^2$ for any two random variables $X$ and $Y$.
\section{Proofs for SBFW}\label{proof:sbfw}
Before proceeding towards the analysis, we provide a Lemma regarding the property of the bias induced by the gradient estimate.
\begin{lemma}\label{old_lemmanew}
	Under Assumption \ref{ass:SBFW}, consider the estimator defined in \eqref{biased_estimate}, then
	
	%	\begin{enumerate}[(a)]
	%		\item \label{parta} 
	(i) define bias $B_t$$:=$$\mbE[h(\x_t,\y_t;\theta_t,\xi_t)]$$-$$\n S(\x_t,\y_{t})$, it holds that we have,
	\begin{align}
	&\norm{B_t}\leq({C_{\x\y}C_{\y}}/{\mu_g})\left(1-({\mu_g}/{L_g})\right)^k,  \label{new1}\\
	& \mbE\|h(\x_t,\y_t;\theta_t,\xi_t)\!-\!\n S(\x_t,\y_{t})\!-\!\!B_t\|^2\leq \sigma_f^2,\label{new2}
	\end{align}
	where $\sigma_f^2$$=$$\sigma_{\x}^2+\frac{3}{\mu_g^2}\left[(\sigma_{\y}^2\!+\!C_{\y}^2)(\sigma_{\x\y}^2\!+\!2C_{\x\y}^2)\!+\!\sigma_{\y}^2C_{\x\y}^2\right]$.
	%\item 
	%\label{partb} 
	
	(ii) For $t\geq 0$, it is possible select $k$ (required to
	approximate the Hessian inverse in \eqref{M}) such that $\|B_t\|\leq \beta_t$ where $\beta_t\leq ct^{a}$ for some constant $c$ and $a>0$.  
	%\end{enumerate}
	
\end{lemma}
%
%\begin{align}
%&\norm{B_t}\leq\frac{C_{\x\y}C_{\y}}{\mu_g}\left(1-\frac{\mu_g}{L_g}\right)^k, \ \ \text{and} \ \ \label{new1}\\
%& \mbE[\|h(\x_t,\y_t;\theta_t,\xi_t)\!-\!\n S(\x_t,\y_{t})\!-\!B_t\|^2]\leq \sigma_f^2.\label{new2}
%\end{align}
\textbf{Proof:} For proof of Lemma \eqref{old_lemmanew}(i) see [Lemma 11, \cite{hong2020two}]. The proof Lemma \eqref{old_lemmanew}(ii) is straight forward. From \eqref{new1} we have $\beta_t=\left(\mathcal{O}(1-\mu_g/L_g)^k\right)$. Now on setting $k=\mathcal{O}(log(t))$ we can get the required condition as $\beta_t\leq ct^{a}$. It shows that with proper selection of $k$, we can make the bias to decay polynomially to zero.

\subsection{Proof of Lemma \ref{Lemma: ySBFW}}\label{proof:Lemma: ySBFW}
Let us consider the term $\mbE_t[\norm{\y_{t}-\y^\star(\x_{t-1})}^2]$ and from the update step 2 of Algorithm \ref{alg:SBFW}, we can write
\begin{align}
&\mbE_t[\norm{\y_{t}-\y^\star(\x_{t-1})}^2]\nonumber\\&=\mbE_t[\norm{\y_{t-1}-\delta_t \n_{\y} g(\x_{t-1},\y_{t-1},\xi_t)-\y^\star(\x_{t-1})}^2].
\end{align}
By expanding the square and taking conditional expectation term inside the inner product terms, we obtain
\begin{align}
&\mbE_t[\norm{\y_{t}-\y^\star(\x_{t-1})}^2]\nonumber
\\&=\mbE_t[\norm{\y_{t-1}-\delta_t \n_{\y} g(\x_{t-1},\y_{t-1},\xi_t)-\y^\star(\x_{t-1})}^2]\nonumber
\\&=\mbE_t[\norm{\y_{t-1}-\y^\star(\x_{t-1})}^2]+\delta_t^2\mbE_t[\| \n_{\y} g(\x_{t-1},\y_{t-1},\xi_t)\|^2]-2\delta_t\mbE_t[\ip{\y_{t-1}\!-\!\y^\star(\x_{t-1}),\n_{\y} g(\x_{t-1},\y_{t-1},\xi_t)}]\label{t1}
\\&\leq \mbE_t[\norm{\y_{t-1}-\y^\star(\x_{t-1})}^2]+\delta_t^2\mbE_t[\| \n_{\y} g(\x_{t-1},\y_{t-1},\xi_t)\|^2]
-2\delta_t\mu_g\mbE_t[\|\y_{t-1}-\y^\star(\x_{t-1})\|^2]\label{t2}\\
&= (1-2\delta_t\mu_g)\mbE_t[\norm{\y_{t-1}-\y^\star(\x_{t-1})}^2]+\delta_t^2\mbE_t[\| \n_{\y} g(\x_{t-1},\y_{t-1},\xi_t)\|^2]\label{y1}
\end{align}
here \eqref{t1} comes from the fact that $\mbE_t[\n_{\y} g(\x_{t-1},\y_{t-1},\xi_t)]=\n_{\y} g(\x_{t-1},\y_{t-1})$, while \eqref{t2} comes from using the strong convexity property of function $g$.
Now let us consider the last term $\mbE_t\| \n_{\y} g(\x_{t-1},\y_{t-1},\xi_t)\|^2$ in the right hand side of \eqref{y1}:
\begin{align}\label{y2}
&\mbE_t[\| \n_{\y} g(\x_{t-1},\y_{t-1},\xi_t)\|^2]\\&=\mbE_t\| \n_{\y} g(\x_{t-1},\!\y_{t-1},\!\xi_t)\!+\!\n_{\y} g(\x_{t-1},\!\y_{t-1})\!\nonumber\\&\quad-\!\n_{\y} g(\x_{t-1},\!\y_{t-1})\|^2\nonumber\\
&\leq 2\mbE_t[\| \n_{\y} g(\x_{t-1},\y_{t-1},\xi_t)-\n_{\y} g(\x_{t-1},\y_{t-1})\|^2]\nonumber\\&\quad+2\|\n_{\y} g(\x_{t-1},\y_{t-1})\|^2,
\end{align}
where we use the inequality $\|\mathbf{a}+\mathbf{b}\|^2\leq 2\|\mathbf{a}\|^2+2\|\mathbf{b}\|^2$. From Assumption \ref{ass:SBFW}(\ref{assSBFW:bounded variance}), we can further upper bound \eqref{y2} as
\begin{align}\label{ohoo}
&\mbE_t[\| \n_{\y} g(\x_{t-1},\y_{t-1},\xi_t)\|^2]\nonumber\\&\leq 2\sigma^2_{g}(1+\norm{\n_{\y}g(\x_{t-1},\y_{t-1})}^2)+2\|\n_{\y} g(\x_{t-1},\y_{t-1})\|^2\nonumber\\
&=2\sigma_g^2+2(1+\sigma_g^2)\|\n_{\y} g(\x_{t-1},\y_{t-1})\|^2\\
&\leq2\sigma_g^2\!+\!2(1\!+\!\sigma_g^2)\|\n_{\y}\! g(\x_{t-1},\y_{t-1})\!-\!\n_{\y} g(\x_{t-1},\!\y^\star(\x_{t-1}))\!\|^2 \nonumber\\
&\leq 2\sigma_g^2+2(1+\sigma_g^2)L_g^2\|\y_{t-1}-\y^\star(\x_{t-1})\|^2
\end{align}
where we used the fact that $\n_{\y} g(\x_{t-1},\y^\star(\x_{t-1}))=0$. Substituting the upper bound in \eqref{ohoo} in  \eqref{y1} we obtain 
\begin{align}\label{y3}
\mbE_t[\norm{\y_{t}-\y^\star(\x_{t-1})}^2]
&\leq [(1-2\delta_t\mu_g)+ 2\delta_t^2(1+\sigma_g^2)L_g^2]\mbE_t[\norm{\y_{t-1}-\y^\star(\x_{t-1})}^2]+2\delta_t^2\sigma_g^2 \nonumber\\
& \leq (1-\delta_t\mu_g)\mbE_t[\norm{\y_{t-1}-\y^\star(\x_{t-1})}^2]+2\delta_t^2\sigma_g^2.
\end{align}
The last inequality in \eqref{y3} is obtained by selecting $\delta_t$ such that $2\delta_t(1+\sigma_g^2)L_g^2\leq \mu_g$. To proceed next, we use Young's inequality to bound the term $\mbE_t[\norm{\y_t-\y^\star(\x_t)}^2]$ in \eqref{y3} as
\begin{align}\label{y4}
\mbE_t[\norm{\y_t-\y^\star(\x_t)}^2]&\leq \left(1+\tfrac{1}{\alpha}\right)\mbE_t[\norm{\y_t-\y^\star(\x_{t-1})}^2]+(1+\alpha)\mbE_t[\norm{\y^\star(\x_t)-\y^\star(\x_{t-1})}^2]
\nonumber\\
&\leq \left(1+\tfrac{1}{\alpha}\right)\mbE_t[\norm{\y_t-\y^\star(\x_{t-1})}^2]+(1+\alpha)\left(\tfrac{C_{\x\y}}{\mu_g}\right)^2\mbE_t\norm{\x_t-\x_{t-1}}^2]\nonumber\\
&\leq \left(1+\tfrac{1}{\alpha}\right)\!\mbE_t[\norm{\y_t\!-\!\y^\star(\x_{t-1})}^2]\!+\!(1\!+\!\alpha)\left(\tfrac{C_{\x\y}}{\mu_g}\right)^2\eta_{t-1}^2D^2
\end{align}
where the second inequality comes from Lemma \ref{Lemma:various_smoothness}(b), and the last inequality comes from the update equation  in \eqref{projected_GD22_4} and the compactness of the domain $\mathcal{X}$. Utilizing \eqref{y3}  into  \eqref{y4}, we get
\begin{align}\label{y5}
&\mbE_t[\norm{\y_t-\y^\star(\x_t)}^2] \nonumber\\&\leq\left(1+\frac{1}{\alpha}\right)(1-\delta_t\mu_g)\mbE_t[\norm{\y_{t-1}-\y^\star(\x_{t-1})}^2] 
+\left(1+\frac{1}{\alpha}\right) 2\delta_t^2\sigma_g^2+(1+\alpha)\left(\frac{C_{\x\y}}{\mu_g}\right)^2\!\!\eta_{t-1}^2D^2.
\end{align}
To proceed next, we substitute $\alpha=\frac{2(1-\delta_t\mu_g)}{\delta_t\mu_g}$ which also implies that $\left(1+\frac{1}{\alpha}\right)(1-\delta_t\mu_g)=1-\frac{\mu_g\delta_t}{2}$:
\begin{align}\label{here3}
\mbE_t[\norm{\y_t-\y^\star(\x_t)}^2]& \leq \left(1-\frac{\delta_t\mu_g}{2}\right)\mbE_t[\norm{\y_{t-1}-\y^\star(\x_{t-1})}^2]\nonumber\\&+\frac{2-\delta_t\mu_g}{\delta_t\mu_g}\left(\frac{C_{\x\y}}{\mu_g}\right)^2\eta_{t-1}^2D^2+\left(1+\frac{1}{\alpha}\right)2\delta_t^2\sigma_g^2\nonumber\\
&\leq  \left(1-\tfrac{\delta_t\mu_g}{2}\right)\mbE_t\norm{\y_{t-1}-\y^\star(\x_{t-1})}^2\!+\!\tfrac{2\eta_{t-1}^2}{\delta_t\mu_g}\left(\tfrac{C_{\x\y}}{\mu_g}\right)^2\!D^2\!+\!4\delta_t^2\sigma_g^2,
\end{align}
where the second inequality comes from the fact that $\frac{2-\delta_t\mu_g}{\delta_t\mu_g}<\frac{2}{\delta_t\mu_g}$ while in the last inequality, we have assumed that $\delta_t$ is chosen such that  $\delta_t\leq \frac{2}{3\mu_g}$ giving $1+\frac{1}{\alpha}\leq 2$. In Corollary \ref{coroSBFW} we will see that our choice of step sizes satisfies these conditions.  
\subsection{Proof of Lemma \ref{Lemma:d:SBFW}} \label{Proof:Lemma:dSBFW}
Starting with update equation \eqref{momentum-track} and employing Lemma \ref{Lemma:general_tracking} we can write 
\begin{align}\label{dSBFW1st}
\mbE_t[\norm{\d_t-\n S(\x_t,\y_{t})-B_t}^2
&\leq(1-\rho_t)^2\mbE_t[\|(\d_{t-1}-\n S(\x_{t-1},\y_{t-1})-B_{t-1})\|] \nonumber\\&\quad+  2 (1-\rho_t)^2\mbE_t[\|h(\x_{t},\y_{t};\theta_t,\xi_t)-h(\x_{t-1},\y_{t-1};\theta_t,\xi_t)\|^2]\nonumber\\
&\quad +2\rho_t^2\mbE_t[\|h(\x_{t},\y_{t};\theta_t,\xi_t)-\n S(\x_t,\y_{t})-B_t)\|^2]\nonumber\\
&\leq(1-\rho_t)^2\mbE_t[\|(\d_{t-1}-\n S(\x_{t-1},\y_{t-1})-B_{t-1})\|] \nonumber\\&\quad+2 \mbE_t[\|h(\x_{t},\y_{t};\theta_t,\xi_t)-h(\x_{t-1},\y_{t-1};\theta_t,\xi_t)\|^2+2\rho_t^2\sigma_f^2,
\end{align}
here the last inequality is obtained using \eqref{new2} and the fact that $(1-\rho_t^2)\leq 1.$ Now we introduce  $h(\x_{t},\y_{t-1};\theta_t,\xi_t)$ and bound the second term of RHS of \eqref{dSBFW1st} as   
\begin{align}\label{derive3}
&\mbE_t[\|h(\x_{t},\y_{t};\theta_t,\xi_t)-h(\x_{t-1},\y_{t-1};\theta_t,\xi_t)\|^2]\nonumber\\
&=\mbE_t\|h(\x_{t},\y_{t};\theta_t,\xi_t)-h(\x_{t-1},\y_{t-1};\theta_t,\xi_t)+h(\x_{t},\y_{t-1};\theta_t,\xi_t)-h(\x_{t},\y_{t-1};\theta_t,\xi_t)\|^2\nonumber\\&\overset{(a)}{\leq} 2\mbE_t[\|h(\x_{t},\y_{t};\theta_t,\xi_t)-h(\x_{t},\y_{t-1};\theta_t,\xi_t)\|^2+2\mbE_t[\|h(\x_{t-1},\y_{t-1};\theta_t,\xi_t)-h(\x_{t},\y_{t-1};\theta_t,\xi_t)\|^2\nonumber\\
&\overset{(b)}{\leq} 2L_k\mbE_t\|\y_{t}-\y_{t-1}\|^2+ 2L_k\mbE_t\|\x_{t-1}-\x_t\|^2\nonumber\\&
\overset{(c)}{\leq} 2L_k\delta_t^2\mbE_t\norm{\n_{\y} g(\x_{t-1},\y_{t-1},\xi_t)}^2+2L_k\eta_{t-1}^2D^2\nonumber\\&
\overset{(d)}{\leq} 2L_k\delta_t^2L_g^2+2L_k\eta_{t-1}^2D^2,
\end{align}
here (a) comes from simple norm property, (b) comes from Lemma \ref{Lemma:borw_result2},
(c) comes from update equation \eqref{inner_update} while (d) comes from Lipschitz continuous Assumption \ref{ass:SBFW}(\ref{assSBFW:Lipschitz_y}) and from the compactness of the set. Using \eqref{derive3} in \eqref{dSBFW1st}, we get the desired expression. 
\subsection{Proof Corollary \ref{coroSBFW}}\label{Proof:corollary:SBFW}   
For the simplicity of analysis we start with writing Lemma \ref{Lemma: ySBFW} for $t=t+1$ and set $\delta_t=\frac{2a_0}{t^q}$ where $a_0=\min\{\frac{1}{3\mu_g},\frac{\mu_g}{2(1+\sigma_g^2)L_g^2}\}$ and $\eta_t=\frac{2}{(t+1)^{\frac{3q}{2}}}$,  which gives
\begin{align}\label{y6}
\mbE_t\norm{\y_{t+1}-\y^\star(\x_{t+1})}^2
&\leq \left(1-\tfrac{2a_0}{(t+1)^q}\right)\mbE_t\norm{\y_{t}-\y^\star(\x_{t})}^2\nonumber\\&\quad+\frac{2}{(t+1)^{3q-q }}\left(\frac{C_{\x\y}}{\mu_g}\right)^2D^2+\frac{16a_0^2}{(t+1)^{2q}}\sigma_g^2\\
&= \left(1-\tfrac{2a_0}{(t+1)^q}\right)\mbE_t[\norm{\y_{t}-\y^\star(\x_{t})}^2]+\tfrac{2(C_{\x\y}/\mu_g)^2D^2+16a_0^2\sigma_g^2}{(t+1)^{2q }}.\nonumber
\end{align}
Note such selection of $\delta_t$ ensures that the conditions  $2\delta_t(1+\sigma_g^2)L_g^2\leq \mu_g$ and $\delta_t\leq \frac{2}{3\mu_g}$ required in Lemma \ref{Lemma: ySBFW} are satisfied. Now taking full expectation and using Lemma \ref{Lemma:gen_seq} we get
\begin{align}\label{final}
\mbE[\norm{\y_{t}-\y^\star(\x_t)}]\leq \frac{b_1}{(t+1)^{q}},
\end{align}
where $b_1=\max\{2^q\norm{\y_1-\y^\star(\x_1)}^2,(2(C_{\x\y}/\mu_g)^2D^2+16a_0^2\sigma_g^2)/(2a_0-1)\}$. Similarly, in Lemma \ref{Lemma:d:SBFW} setting $\delta_t=\frac{2a_0}{(t)^q}$, $\eta_t=\frac{2}{(t+1)^{\frac{3q}{2}}}$ and $\rho_t=\frac{2}{(t)^q}$, we can write

\begin{align}
\mbE_t[\|\d_{t+1}-\n S(\x_{t+1},\y_{t+1})-B_{t+1}\|^2]&\leq \left(1-\frac{2}{(t+1)^q}\right)\mbE_t[\|(\d_{t}-\n S(\x_{t},\y_{t})-B_{t})\|^2]\nonumber\\&+\frac{16L_kL_g^2}{(t+1)^{2q}}+\frac{16L_kD^2}{(t+1)^{3p}}+\frac{8\sigma_f^2}{(t+1)^{2q}}\nonumber\\&\leq \left(1-\frac{2}{(t+1)^q}\right)\mbE_t[\|(\d_{t}-\n S(\x_{t},\y_{t})-B_{t})\|^2]\nonumber\\&+\frac{16L_kL_g^2+16L_kD^2+8\sigma_f^2}{(t+1)^{2q}},
\end{align}
here  the last inequality is obtained using the fact $1/(t+1)^{3q}\leq 1/(t+1)^{2p} $.
Application of Lemma \ref{Lemma:gen_seq} gives
\begin{align}\label{oho2}
\mbE_t[\norm{\d_t-\n S(\x_t,\y_{t})-B_t}^2]\leq \frac{b_2}{(t+2)^{q}},
\end{align}
where $b_2=\max\{2^{q}\norm{\d_1-\n S(\x_1,\y_1)-B_1}^2, 8(2L_kL_g^2+L_kD^2+\sigma_f^2)\}=8(2L_kL_g^2+L_kD^2+\sigma_f^2)$. As, we have initialize $\d_1=h(\x_1,\y_1;\theta_1,\xi_1) $ we can use the bound $$\norm{\d_1-\n S(\x_1,\y_1)-B_1}^2=\norm{h(\x_1,\y_1;\theta_1,\xi_1)-\n S(\x_1,\y_1)-B_1}^2\leq \sigma_f^2$$.
Next, we can bound the term $\mbE\norm{\n Q(\x_t)-\d_t}^2$ as follows
\begin{align}
&\mbE\norm{\n Q(\x_t)-\d_t}^2\\&= \mbE\norm{\n Q(\x_t)-\d_t+B_t+\n S(\x_t,\y_{t})-B_t-\n S(\x_t,\y_{t})}^2 \nonumber\\&\leq 3\mbE\norm{\n Q(\x_t)-\n S(\x_t,\y_{t})}^2+3\norm{B_t}^2+3\mbE\norm{\n S(\x_t,\y_{t})+B_t-\d_t}^2\nonumber\\
&\leq 3\mbE\norm{\y^\star(\x_t)-\y_{t}}^2+3\beta_t^2+\frac{3b_2}{(t+1)^{q}}\nonumber\\&\leq \frac{3b_1}{(t+1)^{q}}+\frac{3b_3}{(t+1)^q}+\frac{3b_2}{(t+1)^{q}}:=\frac{C_1}{(t+1)^{q} },
\end{align}
here second inequality comes from simple norm property, while third inequality is obtained using Lemma \eqref{Lemma:various_smoothness}(a) on the first term, Lemma \ref{old_lemmanew} on the second term and \eqref{oho2} on the third term. The last inequality comes from \eqref{final} and using $\beta_t\leq \frac{C_{\x\y}C_{\y}}{\mu_g(t+1)^q}:=\frac{b_3}{(t+1)^q}$ and the constant $C_1=3(b_1+b_2+b_3)$ is defined as $$C_1=3(\max\{2^q\norm{\y_1-\y^\star(\x_1)}^2,(2(C_{\x\y}/\mu_g)^2D^2+16a_0^2\sigma_g^2)/(2a_0-1)\}+8(2L_kL_g^2+L_kD^2+\sigma_f^2)+\tfrac{C_{\x\y}C_{\y}}{\mu_g}).$$
\section{Proofs for SCFW}
Defining $\tilde{\n} C(\x_t,\y_{t})=\mbE[\n C(\x_t,\y_t,\thh_t,\xi_t)]$, we provide an upper bound on $\mbE\norm{\d_t-\tilde{\n} C(\x_t,\y_{t})}^2$ as follows:
 \begin{lemma}\label{Lemma: f-d}
 	For the iterates $\x_t$ generated by the proposed Algorithm \ref{alg:SCFW}, we have following upper bound
 	\begin{align}
 	\mbE\|\d_t-\tilde{\n} C(\x_t,\y_{t})\|^2&\leq(1-\rho_t)^2\mbE\|\d_{t-1}-\tilde{\n} C(\x_{t-1},\y_{t-1})\|^2 \nonumber\\&+4\rho_t^2M_fM_h+4(1-\rho_t)^2\big[(M_f L_h+3M_h^2L_f)\eta_{t-1}^2 D^2\nonumber\\&+3M_hL_f\delta_t^2\mbE_t\big[\|\y_{t-1}-h(\x_{t-1})\|^2\big]+3\delta_t^2 M_hL_f \sigma^2_h\big],\nonumber
 	\end{align}
 \end{lemma}     
\textbf{Proof:} Setting $\Psi=\tilde{\n} C(\x_t,\y_{t})$ in Lemma \eqref{Lemma:general_tracking}, we can write
\begin{align}\label{main lemma dt}
\mbE_t\|\d_t\!-\!\tilde{\n} C(\x_t,\y_{t})\|^2&=(1\!-\!\rho_t)^2\|\d_{t-1}\!-\!\tilde{\n} C(\x_{t-1},\y_{t-1})\|^2\nonumber\\& \quad+2(1-\rho_t)^2\mbE_t\big[\|\n C_{\phi_{t}}(\x_t,\y_{t})-\n C_{\phi_{t}}(\x_{t-1},\y_{t-1})\|^2\big]\nonumber\\&\quad+2\rho_t^2\mbE_t\big[\|\n C_{\phi_{t}}(\x_t,\y_{t})-\tilde{\n} C(\x_t,\y_{t})\|^2\big],
\end{align} 
We can bound the second term of \eqref{main lemma dt} by taking full expectation and using Assumption \ref{ass:SCFW} as follows
\begin{align}\label{second ka second}
&2\rho_t^2\mbE\|\n C_{\phi_{t}}(\x_t,\y_{t})-\tilde{\n} C(\x_t,\y_{t})\|^2\nonumber\\&=2\rho_t^2\mbE\| \n\hht(\x_t)^\top \n \ft(\y_{t})-\n h(\x_t)^\top \n f(\y_{t})\|^2\nonumber\\&\leq 2\rho_t^2\big(2\mbE\norm{ \n\hht(\x_t)}^2  \mbE\| \n \ft(\y_{t})\|^2+2\mbE\norm{\n h(\x_t)}^2\mbE\norm{ \n f(\y_{t})}^2\big)\leq 4\rho_t^2M_fM_h.
\end{align}
We can bound the first term of \eqref{main lemma dt} by introducing $\n\hht(\x_{t-1})^\top \n \ft(\y_{t})$ as follows
\begin{align}\label{second_term_half}
&\mbE_t[\norm{\n C_{\phi_{t}}(\x_t,\y_{t})-\n C_{\phi_{t}}(\x_{t-1},\y_{t-1})}^2]\nonumber\\
&= \mbE_t[\| \n\hht(\x_t)^\top \n \ft(\y_{t})-\n\hht(\x_{t-1})^\top \n \ft(\y_{t-1}) \|^2]\nonumber\\
&= \mbE_t[\| \n\hht(\x_t)^\top \n \ft(\y_{t})-\n\hht(\x_{t-1})^\top \n \ft(\y_{t}) \nonumber\\&\quad-\n\hht(\x_{t-1})^\top \n \ft(\y_{t-1})+\n\hht(\x_{t-1})^\top \n \ft(\y_{t})\|^2]\nonumber\\
&= \mbE_t[\| \n \ft(\y_{t})^\top(\n\hht(\x_t) -\n\hht(\x_{t-1}))\nonumber\\&\quad+\n\hht(\x_{t-1})^\top( \n \ft(\y_{t})-\n \ft(\y_{t-1}))\|^2]\nonumber\\
&\leq 2\mbE_t[\norm{\n \ft(\y_{t})}^2]\mbE_t[\norm{\n\hht(\x_t) -\n\hht(\x_{t-1})}^2]\nonumber\\&\quad+2\mbE_t[\norm{\n\hht(\x_{t-1})}^2]\mbE_t[\norm{ \n \ft(\y_{t})-\n \ft(\y_{t-1})}^2]\nonumber\\
&\leq 2M_f L_h \mbE_t[\norm{\x_t-\x_{t-1}}^2]+2M_hL_f\mbE_t[\norm{\y_{t}-\y_{t-1}}^2]\nonumber\\
&\leq 2M_f L_h\eta_{t-1}^2 D^2+2M_hL_f\mbE_t[\norm{\y_{t}-\y_{t-1}}^2].
\end{align} We can also bound the second term of \eqref{second_term_half} using the update equation for $y_{t+1}$ from Algorithm:
\begin{align}\label{temp1}
\y_{t}-\y_{t-1}&=-\delta_t \y_{t-1}-(1-\delta_t)\hht(\x_{t-1})+ \hht(\x_t)\nonumber\\&=\delta_t(h(\x_{t-1})-\y_{t-1})+\delta_t(\hht(\x_{t-1})-h(\x_{t-1}))+\hht(\x_t)-\hht(\x_{t-1}).
\end{align}

Taking norm square of \eqref{temp1} we get
\begin{align}\label{temp2}
\mbE_t[\norm{\y_{t}-\y_{t-1}}^2]&\leq 3\delta_t^2\mbE_t\norm{\y_{t-1}\!-\!h(\x_{t-1})}^2\!+\!3\delta_t^2 \sigma^2_h\!+\!3M_h\mbE_t\norm{\x_t\!-\!\x_{t-1}}^2\nonumber\\
&\leq 3\delta_t^2\mbE_t\norm{\y_{t-1}-h(\x_{t-1})}^2+3\delta_t^2 \sigma^2_h+3M_h\eta_{t-1}^2D^2.
\end{align}
Substituting \eqref{temp2} in \eqref{second_term_half} we get
\begin{align}\label{second_term_half2}
&\mbE_t\big[\|\n C_{\phi_{t}}(\x_t,\y_{t+1})-\n C_{\phi_{t}}(\x_{t-1},\y_t)\|^2\big]\\
&\leq 2M_f L_h\eta_{t-1}^2 D^2+2M_hL_f\bigg(3\delta_t^2\mbE_t\big[\|\y_{t-1}-h(\x_{t-1})\|^2\big]+3\delta_t^2 \sigma^2_h+3M_h\eta_{t-1}^2D^2\bigg)\nonumber\\
&=2(M_f L_h+3M_h^2L_f)\eta_{t-1}^2 D^2+6M_hL_f\delta_t^2\mbE_t\big[\|\y_{t-1}-h(\x_{t-1})\|^2\big]+6\delta_t^2 M_hL_f \sigma^2_h.
\end{align}
Thus using bound from \eqref{second ka second} and \eqref{second_term_half2} into  \eqref{main lemma dt} and taking full expectation we obtain
\begin{align}\label{main lemma dt2}
\mbE\|\d_t\!-\!\tilde{\n} C(\x_t,\y_{t})\|^2\!&\leq(1\!-\!\rho_t)^2\mbE\|\d_{t-1}\!-\!\tilde{\n} C(\x_{t-1},\y_{t-1})\|^2 \nonumber\\&\quad+4\rho_t^2M_fM_h+4(1-\rho_t)^2\big[(M_f L_h+3M_h^2L_f)\eta_{t-1}^2 D^2\nonumber\\&\quad+3M_hL_f\delta_t^2\mbE_t\|\y_{t-1}-h(\x_{t-1})\|^2+3\delta_t^2 M_hL_f \sigma^2_h\big].
\end{align}
\subsection{Proof of Corollary \ref{coro1}}\label{proof:coro1}
Starting with Lemma \ref{Lemma:ySCFW} and using the fact that $(1-\delta_t)^2\leq (1-\delta_t) \leq 1$ and substituting $\delta_t=\frac{2}{t^{p}}$ and $\eta_t\leq\frac{2}{(t+1)^{p}}$, we can write
\begin{align}\label{ytemp2}
&\mbE\norm{\y_{t+1}-h(\x_{t+1})}^2\leq \left(1-\tfrac{2}{(t+1)^{p}}\right)\mbE\norm{(\y_{t}-h(\x_{t}))}^2+\tfrac{8\sigma_h^2}{(t+1)^{2p}}+\tfrac{8M_hD^2}{(t+1)^{2p}}.
\end{align}	
Application of Lemma  \ref{Lemma:gen_seq} in \eqref{ytemp2} gives
\begin{align}\label{bound_y-g2}
\mbE\norm{\y_{t}-h(\x_{t})}^2\leq \frac{a_1}{(t+1)^{p}},
\end{align}
where $a_1=\max\{2^{p}\mbE\norm{\y_1-h(\x_1)}^2, 8(\sigma_h^2+M_hD^2)\}$. As we have initialized $\y_1=g_{\xi_1}(\x_1)$ we can bound $\mbE\norm{\y_1-h(\x_1)}^2=\mbE\norm{g_{\xi_1}(\x_1)-h(\x_1)}^2\leq \sigma_h^2$. Further, using the fact that $0<p\leq 1$, we can simplify $a_1$ as  $a_1=8(\sigma_h^2+M_hD^2)$. Similarly,
from Lemma \ref{Lemma: f-d} using the fact $(1-\rho_t)^2\leq 1$, we can write
\begin{align}\label{main lemma dt32}
&\mbE\|\d_{t+1}-\tilde{\n} C(\x_{t+1},\y_{t+1})\|^2\\&\leq(1-\rho_{t+1})^2\mbE\|\d_{t}-\tilde{\n} C(\x_{t},\y_{t})\|^2 +4\rho_{t+1}^2M_fM_h\nonumber\\&\quad+4\big[(M_f L_h+3M_h^2L_f)\eta_{t}^2 D^2+3M_hL_f\delta_{t+1}^2\mbE_t\big[\|\y_{t}-h(\x_{t})\|^2\big]+3\delta_{t+1}^2 M_hL_f \sigma^2_h\big].
\end{align}
Substituting $\delta_t=\frac{2}{t^{p}}$, $\rho_t=\frac{2}{t^{p}}$ and $\eta_t\leq\frac{2}{(t+1)^{p}}$, and using bound \eqref{bound_y-g2} we get
\begin{align}\label{main lemma dt42}
&\mbE\|\d_{t+1}-\tilde{\n} C(\x_{t+1},\y_{t+1})\|^2\\&\leq\left(1-\frac{2}{(t+1)^{p}}\right)\mbE\|\d_{t}-\tilde{\n} C(\x_{t},\y_{t})\|^2 \nonumber\\&\quad+16\bigg[\frac{M_fM_h}{(t+1)^{2p}}+\frac{(M_f L_h+3M_h^2L_f)D^2}{(t+1)^{2p}}+\left(\frac{3M_hL_f}{(t+1)^{2p}}\right)\left(\frac{a_1}{(t+1)^{p}}\right)+\frac{3 M_hL_f \sigma^2_h}{(t+1)^{2p}}\bigg].
\end{align}
Using Lemma  \ref{Lemma:gen_seq} and the fact that $\tfrac{1}{(t+1)^{3p}}\leq \tfrac{1}{(t+1)^{2p}}$ gives
\begin{align}\label{main lemma dt52}
&\mbE\|\d_t-\tilde{\n} C(\x_t,\y_{t})\|^2\leq\frac{a_2}{(t+1)^{p}},
\end{align}
where $a_2:=\max\{2^{p}\mbE\|\d_{1}-\tilde{\n} C(\x_{1},\y_{1})\|^2,16(M_fM_h+(M_f L_h+3M_h^2L_f)D^2+3M_hL_fa_1+3 M_hL_f \sigma^2_h)\}$. As we have initialized $\d_{1}=\n C_{\phi_{1}}(\x_1,\y_{1})$, we can bound $\mbE\|\d_{1}-\tilde{\n} C(\x_{1},\y_{1})\|^2=\mbE\|\n C_{\phi_{1}}(\x_1,\y_{1})-\tilde{\n} C(\x_{1},\y_{1})\|^2\leq 2M_hM_f$ (see \eqref{second ka second}). Also as $0<p\leq1$, we can further simplify $a_2$ to get $a_2=16(M_fM_h+(M_f L_h+3M_h^2L_f)D^2+3M_hL_fa_1+3 M_hL_f \sigma^2_h)$.

Now, introducing $\tilde{\n} C(\x_t,\y_{t})$ can write
\begin{align}\label{d-gradF2}
&\mbE\|\n C(\x_t)-\d_t\|^2\nonumber\\&=\mbE\|\n C(\x_t)-\d_t-\tilde{\n} C(\x_t,\y_{t})+\tilde{\n} C(\x_t,\y_{t})\|^2 \nonumber\\& \leq 2\mbE\|\d_t-\tilde{\n} C(\x_t,\y_{t})\|^2+2\mbE\|\tilde{\n} C(\x_t,\y_{t})-\n C(\x_t)\|^2\nonumber\\& = 2\mbE\|\d_t-\tilde{\n} C(\x_t,\y_{t})\|^2+2\mbE\norm{\n h(\x_t)^\top \n f(\y_{t})-\n h(\x_t)^\top\n f(h(\x_t))}^2\nonumber\\& \leq 2\mbE\|\d_t-\tilde{\n} C(\x_t,\y_{t})\|^2+2M_hL_f\mbE\|\y_{t}-h(\x_t)\|^2\nonumber\\& \leq \frac{2a_2}{(t+1)^{p}}+\frac{2M_hL_fa_1}{(t+1)^{p}}:=\frac{A_1}{(t+1)^{p}},
\end{align}
here the first inequality comes from simple norm property while in the next expression to it we simply substitute the values of gradients. In the second inequality, we used the Assumption \ref{ass:SCFW} (\ref{ass:smoothness}-\ref{ass:bounded_moment}) while the last inequality is obtained using bounds from \ref{bound_y-g2} and  \ref{main lemma dt52} and defined $A_1:=32[M_h(M_f+28L_f\sigma_h^2)+(M_f L_h+28M_h^2L_f)D^2]$.     
\subsection{Proof of Theorem \ref{theorem1}}\label{proof:theorem1}
\noindent\subsubsection{Proof of Statement (i) (convex case)}
Starting with the smoothness of $C$ with parameter $L_F=M_h^2L_f + M_fL_h$ and proceeding the same way as \eqref{smoothQ}-\eqref{Lemma for Q final eq2} we can write
\begin{align}\label{Lemma for F final eq2}
\mbE[C(\x_{t+1})-C(\x^\star)]
&\leq (1-\eta_{t})\mbE[C(\x_t)-C(\x^\star)]\nonumber\\&\quad+ \eta_{t} D\sqrt{\mbE\norm{\n C(\x_t)-\d_t}^2}+\frac{L_F\eta_{t}^2D^2}{2}.
\end{align}
Setting $\eta_t=\frac{2}{t+1}$ and using Corollary \ref{coro1} with $p=1$, we can bound the second term of \eqref{Lemma for F final eq2} as
\begin{align}\label{sqrt d-f}
\eta_tD\sqrt{\mbE\norm{\n C(\x_t)-\d_t}^2} \leq \frac{2D\sqrt{A_1}}{(t+1)^{3/2}}.
\end{align}
Now using bound from \eqref{sqrt d-f} in \eqref{Lemma for F final eq2} we obtain
\begin{align}\label{Lemma for F final eq3}
\mbE[C(\x_{t+1})-C(\x^\star)]&
\leq \left(1-\frac{2}{t+1}\right)\mbE[C(\x_t)-C(\x^\star)]+ \frac{2 D\sqrt{A_1}}{(t+1)^{\frac{3}{2}}}+\frac{2L_FD^2}{(t+1)^{2}}.
\end{align}
Multiplying both side by $t(t+1)$ we can write
\begin{align}\label{Fcc }
&t(t+1)\mbE[C(\x_{t+1})-C(\x^\star)]\\&
\leq t(t-1)\mbE[C(\x_t)-C(\x^\star)]+\frac{ 2t D\sqrt{A_1}}{(t+1)^{\frac{1}{2}}}+\frac{2tL_FD^2}{t+1}\nonumber\\
&\leq t(t-1)\mbE[C(\x_t)-C(\x^\star)]+ 2 D\sqrt{A_1}(t+1)^{\frac{1}{2}}+2L_FD^2,\nonumber
\end{align}
here in last the inequality we used the fact that $t<t+1$. Summing for $t=1,2,\cdots,T$ and rearranging we get
\begin{align}\label{Fcc1}
&\mbE[C(\x_{T+1})-C(\x^\star)]\nonumber\\
&\leq \frac{1}{T(T+1)}\left(2 D\sqrt{A_1}\sum_{t=1}^{T}(t+1)^{\frac{1}{2}}+2L_FD^2T\right)\nonumber\\
&\leq \frac{1}{T(T+1)}\left(\frac{4}{3} D\sqrt{A_1}(T+1)^{\frac{3}{2}}+2L_FD^2T\right)\nonumber\\
&\leq \frac{8 D\sqrt{A_1}}{3(T+1)^{\frac{1}{2}}}+\frac{2L_FD^2}{(T+1)},
\end{align}
here in the second inequality we use the fact that $\sum_{t=1}^{T}(t+1)^{1/2}\leq \frac{2}{3}(T+1)^{3/2}$ and the last inequality is obtained using $(T+1)=T(1+\frac{1}{T})\leq 2T $. 

%$16(M_fM_h+(M_f L_h+3M_h^2L_f)D^2+24M_hL_f(\sigma_h^2+M_hD^2)+3 M_hL_f \sigma^2_h)+8M_hL_f(\sigma_h^2+M_hD^2)\\
%8[2M_fM_h+2(M_f L_h+3M_h^2L_f)D^2+48M_hL_f(\sigma_h^2+M_hD^2)+6 M_hL_f \sigma^2_h+M_hL_f(\sigma_h^2+M_hD^2)]\\
%8[2M_fM_h+2(M_f L_h+3M_h^2L_f)D^2+48M_hL_f\sigma_h^2+48M_hL_fM_hD^2+6 M_hL_f \sigma^2_h+M_hL_f\sigma_h^2+M_h^2L_fD^2]\\
%8[2M_fM_h+(2M_f L_h+55M_h^2L_f)D^2+55M_hL_f\sigma_h^2]\\
%16[M_fM_h+(M_f L_h+28M_h^2L_f)D^2+28M_hL_f\sigma_h^2]\\
%16[M_h(M_f+28L_f\sigma_h^2)+(M_f L_h+28M_h^2L_f)D^2]
%$
%
\subsubsection{Proof of Statement (ii) (Non-convex case)}
Again starting with  the smoothness assumption of $F$ and proceeding the same way as \eqref{smooth non convexQ}-\eqref{sumQ}, we can write
\begin{align}\label{sum}
&\sum_{t=1}^{T}\eta_{t}\mbE[\mathcal{G}(\x_t)]
\leq C(\x_1)-C(\x^{\star})+2D\sum_{t=1}^{T}\eta_{t}\mbE\norm{\d_t-\n C(\x_t)}+\frac{L_FD^2}{2}\sum_{t=1}^{T}\eta_{t}^2.
\end{align}
Using bound from  Corollary \ref{coro1} for $p=2/3$, and  Jensen's equality we can write
\begin{align}\label{jen}
\mbE\norm{\d_t-\n C(\x_t)}\leq\sqrt{\mbE\norm{\d_t-\n C(\x_t)}^2}\leq \tfrac{\sqrt{A_1}}{(t+1)^{\frac{1}{3}}}.
\end{align}
%Note that our selection of $\eta_t=\eta=\frac{2}{(T+2)^{2/3}}$ satisfies the condition $\eta_t\leq\frac{2}{(t+2)^{2/3}}$ required in Lemma \ref{Lemma: ySBFW} and Lemma \ref{Lemma: f-d}. 
Now,
using \eqref{jen}, setting $\eta_t=\frac{2}{(T+1)^{2/3}}$ in \eqref{sum} and rearranging we can write
\begin{align}
&\mbE[\mathcal{G}(\hat{\x})]\leq \frac{1}{T}\sum_{t=1}^{T}\mbE[\mathcal{G}(\x_t)]
\nonumber\\&\leq\frac{C(\x_1)-C(\x^\star)}{2T(T+1)^{-2/3}}+\sum_{t=1}^{T}\frac{2D\sqrt{A_1}}{T(t+1)^{1/3}}+\sum_{t=1}^{T}\frac{L_FD^2}{T(T+1)^{2/3}}\nonumber\\
&\leq \frac{C(\x_1)-C(\x^\star)}{2T(T+1)^{-2/3}}+\frac{3 D\sqrt{A_1}}{T(T+1)^{-2/3}}+\frac{L_FD^2}{(T+1)^{2/3}}\nonumber\\&= \frac{(T\!+\!1)(C(\x_1)\!-\!C(\x^\star))}{2T(T+1)^{1/3}}\!+\!\frac{3 D\sqrt{A_1}(T\!+\!1)}{T(T+1)^{1/3}}\!+\!\frac{L_FD^2}{(T+1)^{2/3}}\nonumber\\
&\leq \frac{C(\x_1)-C(\x^\star)}{(T+1)^{1/3}}+\frac{6 D\sqrt{A_1}}{(T+1)^{1/3}}+\frac{L_FD^2}{(T+1)^{2/3}}\nonumber\\
&:=\mathcal{O}((T+1)^{-1/3}),
\end{align}
here in the second inequality we use the fact that $\sum_{t=1}^{T}(t+1)^{-1/3}\leq \frac{3}{2}(T+1)^{2/3}$ and the last inequality is obtained using $(T+1)=T(1+\frac{1}{T})\leq 2T $.

   \end{document}